 \UseRawInputEncoding
\documentclass[12pt]{amsart}
\usepackage{longtable}
\usepackage{supertabular}
\usepackage[dvips]{graphicx}
\usepackage{graphicx}
\usepackage[colorlinks=true]{hyperref}     
\hypersetup{linkcolor=black, citecolor=black}
\newtheorem{thm}{Theorem}[section]

\newtheorem{prop}[thm]{Proposition}
\newtheorem{rmk}[thm]{Remark}

\theoremstyle{definition}

\theoremstyle{remark}

\numberwithin{equation}{section}

\newcommand{\R}{\mathbb R}
\def\disp{\displaystyle}

\def\argmin{\mathop{{\rm argmin}}\nolimits}

\def\spt{\mathop{{\rm spt}}\nolimits}

\def\h{{\mathcal H}^1}

\def\res{\hbox{\vrule height8pt\vrule height0.4pt width8pt\kern1pt}}

\def\sbh {S\kern-2pt B\kern-2pt H}  
\def\sbv {S\kern-2pt B\kern-1pt V}  \def\sbvo{\sbv(\Omega)}
             
\def\gsbv {G\kern-1pt S\kern-1.5pt BV}        
\def\xx{{\bf x}}

\def\ro{\varrho}

\def\O{\Omega}

\begin{document}
\parindent=0pt
\thispagestyle{empty}
\title[ELASTIC-BRITTLE REINFORCEMENT OF FLEXURAL STRUCTURES]{Elastic-brittle reinforcement\\ of flexural
structures}

\author[MPT]{Francesco Maddalena$^*$, Danilo Percivale$^{**}$
\& Franco Tomarelli$^{***}$}
\address{$^*$ Politecnico di Bari, Dipartimento di Meccanica, Matematica e Management, Italy}%
\address{$^{**}$ Universit\`a di Genova, Dipartimento della Produzione
Termoenergetica e Modelli Matematici, Italy}%
\address{\noindent$^{***}$ Politecnico di Milano, Dipartimento di Matematica, 
Italy  \vskip0.5cm }

\email{f.maddalena@poliba.it}
\email{percivale@diptem.unige.it}
\email{franco.tomarelli@polimi.it}%
\date{\today}

\thanks{This work is supported by PRIN (M.I.U.R)
\textit{``Variational methods for stationary and evolution problems with singularities and interfaces'' 2017}
and Progetto G.N.A.M.P.A. (INDAM) \textit{``Problemi asintotici e fenomeni singolari in meccanica dei continui (2020)''} .}
\maketitle
\rightline{\textsl{Dedicated to Claudio Baiocchi with sore heart for the loss of a teacher and friend}}
\vskip0.6cm
 {\bf Abstract.}
This note provides a variational description of the mechanical effects of
flexural stiffening of a 2D plate glued to
an elastic-brittle or an elastic-plastic reinforcement.
The reinforcement is assumed to be linear elastic outside
possible free plastic yield lines or free crack. Explicit Euler equations and a compliance identity are
shown for the reinforcement of a 1D beam.

\vskip1cm

\tableofcontents

\section{Introduction}\label{intro}
\noindent
The main theme moving the present paper relies in studying the mechanisms ruling stress transfer between  material structures affected by strongly different constitutive properties.
We consider some functionals of the kind 
\begin{eqnarray}
\displaystyle
  &&
\ \ \Lambda(K,v_r,v_p) \ =\ \mathcal{H}^1\!\left( K \cap \overline \Omega\right)\!
  \,+ \,\sigma\!\!  \int_K \! |\,[Dv_r]\,|\,d\mathcal H^1 +
   \!\!\int_{\Omega\setminus K}\! \!\!\big( \eta\,|D^2 v_r|^2 -f_rv_r\big) \,d\xx +               \\
  \nonumber
  \displaystyle
  && \qquad\qquad\qquad\quad
               +\,\mu \int_{\Omega} |v_r-v_p|^2 \, d\xx           
             +  \!\int_{\Omega} \!\left(\gamma\,|D^2 v_p|^2
              \! -  f_p \, v_p \right) \, d\xx
\end{eqnarray}
dependent on competing triplets $(K,v_r,v_p)$.
Here $\Omega\subset\R^2$ 
is given together with the loads $f_r, f_p$, the nonnegative constitutive parameters $\sigma,\,\eta,,\,\mu,\,\gamma$ and a suitable Dirichlet boundary condition, shared by both displacements $v_r$ and $v_p$. $\mathcal H^1$ denotes the 1-dimensional Hausdorff measure. Concerning competing triplets $(K,v_r,v_p)$, 
the set $K$ is  
closed, the function $v_r$ is smooth outside $K$ while $v_p$ is smooth on $\Omega$.
\\
Structural reinforcements are extensively employed in manufactured engineering systems, ranging from the traditional field of composite structures to the more recent applications in microelectronic devices and nano-reinforced composites (\cite{PHUONG nano}).
Indeed, external bonding of plates is a method of strenghtening which involves an additional adhering reinforcement to a structural element (\cite{PPG}). The adhesive is needed to transfer the stresses between the two  elements. This technique is aimed to reduce deflection, hence confine crack and plastic yielding location in order to increase load carrying capacity (see e.g. \cite{CCT}, \cite{HJP}, \cite{KN}, \cite{PSC}, \cite{ZNF}), as well as predicting the behavior of paint coating layers (\cite{PHUONG nano}, \cite{co-de-gecko}).\\

In a  typical
 reinforcement system, represented by a  brittle structure (\cite{DMI}) bonded  to a more compliant substrate,  cracking and debonding instabilities (delamination)
of the brittle element  may appear  under the action of external data
that may be ruled by external loading, temperature change or even residual thermal stresses.
The occurrence of  plastic yielding, cracking and loss of adhesion (or delamination) constitute the main failure modes of reinforcements, so
 a great deal of work has been done over past decades
to apply fracture mechanics in description of behavior and the influence of cracks nucleated on or near an
interface between two dissimilar materials and a number of papers have  been published
on the problem (see, for instance, \cite{CYS}, \cite{HS}, \cite{KN}).
The crucial questions, when studying the possible failure of bonded structures, rely in understanding  crack or yielding nucleation and crack propagation in presence of debonding or delamination of the constituent materials.
This issue, precisely the role played by bonding layer  in the formation of singularities,  has not been yet investigated  in spite of its influence in the mechanical behavior of such structural systems.
Therefore, in our opinion, an appropriate description of such problems should incorporate all these \textsl{strongly nonlinear} effects in a mathematical theory which is able to detect qualitative and quantitative features of the underlying physics (\cite{MPT5},\cite{MPT6}).\\
To this aim we propose in the present paper a variational approach in which debonding and 
possible singular states arise as minimizers of suitable energy functionals.\\
In the previous works  \cite{MP}, \cite{MPT1}, \cite{MPT2}, \cite{MPT3}, 
we have studied the adhesion interaction of linear and nonlinear elastic structures
by focusing on the influence of different constitutive choices for the adhesive material, while in \cite{MPT3} we have investigated the occurrence of global collapse and the interplay of cracking and debonding for a couple of plane elastic
sheets.\\
Here reinforcements are modelled in the framework of Kirchhoff-Love theory,  while the addition of a fracture energy term according to the Griffith theory allows to capture crack formation.  This program is achieved
by exploiting  the techniques developed in the study of second order variational problems
with free  discontinuity
(\cite{WPLATE}
-\cite{SBZ},
\cite{BZEE}-\cite{SBZDIR},
\cite{PT3},\cite{PT1}).
\vskip0.2cm
Precisely we aim to describe
the mechanical effects of bending a 2-dimensional 
stiff substrate fixed at the boundary, shortly called
clamped plate from now on, which is glued to
an elastic-plastic-brittle reinforcement. The reinforcement is assumed to be linear elastic outside
possible free plastic yield lines or free 1D crack.
The adhesive interaction between the structures
is modelled through an energetic contribution whose density is the 
square of the modulus of difference of the displacement;
concerning this assumption we recall that, 
if one assumes 
the $q$-power of the modulus of difference, with
$0<q<1$, then the only stable configurations of the system are the completely detached or completely glued ones,
even in the case of a flat plate as it was proved in Section 2.2 of \cite{MP}.

The phenomenon is modeled as a variational problem which allows
free discontinuity and free gradient discontinuity for the reinforcement.
We assume that
both configurations of the plate and its reinforcement are described as graphs referred to the coordinates in the horizontal plane, undergo vertical displacements and are subject to Dirichlet boundary conditions.
We describe the details in these three main cases.
\begin{enumerate}
  \item \textbf{hard-device reinforcement:} the structure consists in the glueing of a plate with a reinforcement; the plate undergoes a prescribed configuration
    (described by a given displacement) and consequently acts on the reinforcement through the adhesive layer; the reinforcement behaves as a piece-wise Kirchhoff-Love plate since it can develop lower-dimensional singularities of two kinds, plastic yielding (free gradient-discontinuity) and/or crack (free discontinuity).
  \item \textbf{strengthening reinforcement:} the structure consists in the glueing of two objects, which are still labeled as plate and reinforcement; the plate behaves like a Kirchhoff-Love plate, whose unstressed flat reference configuration is horizontal, under the action
  of a given transverse (vertical) load $f$, while the plate displacement acts on the reinforcement through the adhesive layer; the reinforcement behaves as a piece-wise Kirchhoff-Love plate since it can develop lower-dimensional singularities of two kinds, plastic yielding and/or crack.
We denote the admissible vertical displacement of the reinforcement by $v_r$ and
the admissible vertical displacement of the plate   by $v_p$.
\item
 \textbf{elastic-plastic reinforcement:} 
 as like as in case (2), but 
 without crack and with a refinement of the yielding energy
 along the a priori. unknown plastic yield lines.
\end{enumerate}
As far as we know in structural mechanics literature there are few studies (\cite{CLT}) of the interplay of plastic-yielding or fracture with bending bulk energy: here we aim to study this coupling in terms of integral functionals with free discontinuity and free gradient discontinuity by methods introduced in calculus of variations 
(see \cite{AFP},\cite{DCDSA},\cite{JMPA},\cite{CLTadvcv},\cite{SICI}).
\\
In Section \ref{1D}
we discuss the analogous one-dimensional case, say elastic-brittle reinforcement of flexural beams (Theorems \ref{Thm 3.1},\,\ref{Thm 3.2},\,\ref{Thm 3.3}). 
In Sections \ref{Statement}, \ref{hard-device Section} and \ref{strengthening Section}
we exhamine details of the clamped plate reinforcement
(hard-device and strengthening)
showing the existence of energy minimizing solutions (Theorems \ref{Theorem 1},\,\ref{Theorem 2},\,\ref{Theorem 3}). 
\\
In the 1D case (beam reiforcement)
we provide explicit Euler equations (Proposition \ref{EulerThm}), transmission conditions at free-discontinuity set
and compliance identities fulfilled by minimizers\,(Proposition \ref{ComplianceThm}). 
Lack of convexity in these functionals may lead to non uniqueness of minimizers (\cite{BoT1},\cite{BoT2}). However we show uniqueness and hence smoothness in case of small loads: see Theorem \ref{Thm reg min for F1}, Remark \ref{Rmk reg min for E1}. A more detailed analysis of non uniqueness phenomenon is postponed in a forthcoming article.
\\
In present paper we omit the consideration of a unilateral constraint forcing non-interpenetration the plate and its reinforcement: such constraint leads to technical difficulties and substantial problems for existence of strong solutions in the 2D case of the plate, since
Euler equations are replaced by variational inequalities (\cite{BC}); this  issue is postponed to a subsequent paper (see Remark \ref{rmk 5 omega}).
Here we only start the analysis of solutions for the 1D cases of the beam reinforced by a hard device, obtaining a variational inequality coupled with free discontinuity (Propositions \ref{Thm 3.5} and \ref{EulerThmConstr}),
and of the strengthening reinforced beam, obtaining quasi-variational inequalities coupled with free discontinuity (Propositions \ref{Thm reinf1 with unilat constr} and \ref{EulerThmConstrF1}).
\section{Statement of the problem and main results}\label{Statement}
Assume
\vskip-0.5cm
\begin{equation}\label{Omega}
\vspace{-0.3cm}
    \Omega\subset\!\subset\Omega_p\subset\!\subset\mathbb{R}^2
    \qquad \hbox{are bounded and connected $C^2$ open sets} 
\end{equation}
\begin{equation}\label{w}
\vspace{-0.1cm}
      w\in  C^2(\overline{\Omega_p}) \,,
\end{equation}
where
$\Omega_p$ represents the horizontal reference configuration of the plate
and $\Omega$ represents the horizontal reference configuration of the reinforcement and
$w:\Omega_p\mapsto \mathbb{R}$
prescribes the Dirichlet datum of the clamped plate.
\\
In the case of hard-device reinforcement $v:\Omega_p\mapsto \mathbb{R}$ denotes the generic admissible vertical displacement of the reinforcement,
while the vertical displacement $w$ of the plate is prescribed:
since the reinforcing structure has to accomplish a prescribed configuration w, this amounts to deal with a pre-strained  state of the material competing with other energetic terms.
In the case of strengthening reinforcement, $v_r:\Omega_p\mapsto \mathbb{R}$ denotes the generic admissible vertical displacement of the reinforcement,
while the generic admissible vertical displacement of the plate is denoted by $v_p$ where $v_p:\Omega_p\mapsto \mathbb{R}$, where $v_p$ is subject to a Dirichlet-type boundary condition, prescribed on $\Omega_p\setminus\Omega$.
We face the three cases mentioned in the Introduction, by studying the minimization
of suitable energy functionals:
\begin{itemize}
  \item the energy $E$ associated to a hard-device reinforcement, which is dependent on pairs $(K,v)$
where $K$ denotes the damaged region of the reinforcement and $v$ its transversal displacement;
  \item the energy $F$ associated to a strengthening reinforcement, which is dependent on triplets $(K, v_r,v_p)$
where $K$ still denotes the damaged region of the reinforcement (free discontinuity and free gradient-discontinuity) while $v_r$ and $v_p$ denote
respectively the transversal displacements of reinforcement and plate;
 \item the energy $G$ associated to an elastic-plastic reinforcement, which is dependent on triplets $(K, v_r,v_p)$ as above, but $v_r$ may undergo only free gradient-discontinuity on $K$.
\end{itemize}
We state some results related to minimization of these energies;
their proofs are postponed in Sections \ref{hard-device Section} and \ref{strengthening Section}. All functions under exam are real-valued.
\begin{thm}\label{Theorem 1}
\textbf{(hard-device reinforcement)} Assume \eqref{Omega},\eqref{w} and
\begin{equation}\label{eta mu}
    \eta>0\,,\qquad \mu>0\,,
     \qquad f\in L^{4}(\Omega)\,,
\end{equation}
then
there exists a pair $(Z,u)$
minimizing
\begin{equation}\label{E}
    E(K,v)\, =\ \mathcal{H}^1\left( K \cap \overline \Omega\right)\, +
                \int_{\Omega\setminus K}\!\! \Big(\,\eta\, |D^2 v|^2 - f v\,\Big)\, d\xx \,+\,
               \mu\int_{\Omega} |v-w|^2 \, d\xx
   \end{equation}
over \textbf{essential admissible pairs} $(K,v)$, say pairs s.t.
\begin{equation}\label{admissible K v}
\left\{
\begin{array}{l}
    \displaystyle K \ \hbox{is the smallest closed subset of } \mathbb{R}^2 \ \hbox{s.t.} \\
    \displaystyle v\in C^2(\Omega_p\setminus K),
    \ \qquad v\equiv w \ \hbox{a.e in }
    \Omega_p\setminus \overline{\Omega} \,.
\end{array}\right.\end{equation}
Moreover $Z\cap\Omega_p=Z\cap\overline\Omega$ is an $(\mathcal{H}^1,1)$ rectifiable set and
$E(Z,u)<+\infty\,.$
\end{thm}
Here and in the sequel we denote by
$\h$ the $1$-dimensional Hausdorff measure. \\
If $(Z,u)\in\argmin E$, say $(Z,u)$ is an optimal pair among the ones fulfilling
\eqref{admissible K v}, then $Z$
represents the damaged zone of the reinforcement $\Omega$, say the 1D set
where either plastic yielding or fracture occur, and $u$ is the related transverse
displacement of the reinforcement.
%
\begin{thm}\label{Theorem 2}
\textbf{(strengthening reinforcement)} - Assume \eqref{Omega},\eqref{w},\eqref{eta mu} and
\begin{equation}\label{f}
    f_r\in L^4(\Omega),\quad  f_p\in L^2(\Omega)\,,\quad \gamma>0\,.
    \vspace{-0.1cm}
\end{equation}
\vskip0.1cm
Then there exists a triplet $(Z,u_r,u_p):=(Z,U)$
minimizing
\begin{equation}\label{F}
\begin{array}{ll} \vspace{0.2cm}
     F(K,v_r,v_p) \, \!\! &=\, F(K,V) \,:=    
    \vspace{0.1cm}
    \\ 
   & \displaystyle
   \mathcal{H}^1\!\left( K \cap \overline \Omega\right) +
               \eta \int_{\Omega\setminus K}\! 
               \big(|D^2 v_r|^2 -f_rv_r\big) d\xx +
        \vspace{0.1cm}      
         \\ 
              & \displaystyle
              \qquad \qquad +
               \mu\int_{\Omega}\! |v_r-v_p|^2  d\xx\,+
               \int_{\Omega} \!\left(\gamma\,|D^2 v_p|^2
              \! -  f_p \, v_p \right)  d\xx
   \end{array}
   \end{equation}
over \textbf{essential admissible triplets} $(K,v_r,v_p)$, say triplets s.t.
\begin{equation}\label{admissible K vr vp}
\left\{
\begin{array}{l}
     v_p\in H^2(\Omega_p)\,, \qquad v_p \equiv w \hbox{ a.e in }
    \Omega_p\setminus \overline{\Omega}\,,
     \vspace{0.1cm}
     \\
    \displaystyle K \ \hbox{is the smallest closed subset of } \mathbb{R}^2 \ \hbox{s.t.} \\
    \displaystyle v_r\in C^2(\Omega_p\setminus K),
    \ \qquad v_r\equiv w \ \hbox{a.e in }
    \Omega_p\setminus \overline{\Omega} \,.
\end{array}\right.\end{equation}
%
Moreover $Z\cap\Omega_p=Z\cap\overline\Omega$ is an $(\mathcal{H}^1,1)$ rectifiable set and
$F(Z,u_r,u_p)<+\infty\,.$
\end{thm}
If $(Z,U)=(Z,u_r,u_p)$ is an optimal triplet among the ones fulfilling
\eqref{admissible K vr vp}, say it is an essential admissible pair
$(Z,U)\in\argmin F$, then $Z$
represents the damaged zone of the reinforcement $\Omega$, say the 1D set
where either plastic yielding or fracture occur, and $u_r,\ u_p$ respectively are the related
displacement of the reinforcement and the plate.
\begin{thm}\label{Theorem 3}
\textbf{(elastic-plastic reinforcement of flexural plate)} - \\ Assume \eqref{Omega},\eqref{w},\eqref{eta mu}, $\sigma\!>\!0,\
f_r\!\in\! L^s(\Omega),\,s\!>\!2$, and $f_p\!\in\! L^2(\Omega)$.
\\ Then
there is a triplet $(Z,u_r,u_p)=(Z,U)$ minimizing
\begin{eqnarray}\label{G}
    \displaystyle
    && G(K,v_r,v_p)\ =\ 
    G(K,V)\ :=\ 
    \\
    \nonumber 
    && 
    \ =\ \mathcal{H}^1\left( K \cap \overline \Omega\right) \,+\,
     \sigma \int_{K\cap \overline\Omega}|\,[Dv_r]\,|\,d\mathcal H^1
               \,+\,
                \int_{\Omega\setminus K}\!\! \big(\eta |D^2 v_r|^2 - f_r v_r\big)\, d\xx \,+\,
\\  
&&
 \nonumber
  \qquad  \,+\,
                \mu\int_{\Omega} |v_r-v_p|^2 \, d\xx
                +\,
                \int_{\Omega}\! \big(\gamma |D^2 v_p|^2 - f_p v_p\big)\, d\xx
   \end{eqnarray}
over \textbf{elastic-plastic essential admissible triplets} $(K,v_r,v_p)$, say triplets s.t.
\begin{equation}\label{admissible K vr vp G}
\left\{
\begin{array}{l}
    \displaystyle K \ \hbox{is the smallest closed subset of } \mathbb{R}^2 \ \hbox{s.t.} \\
    \displaystyle v_r\in C^0(\Omega_p)\cap C^2(\Omega_p\setminus K),
    \ \qquad v_r= v_p\equiv w \ \hbox{a.e in }
    \Omega_p\!\setminus\! \overline{\Omega} \,.
\end{array}\right.\end{equation}
Moreover $Z\cap\Omega_p=Z\cap\overline\Omega$ is an $(\mathcal{H}^1,1)$ rectifiable set and
$G(Z,u_r,u_p)<+\infty\,.$
\end{thm}
Theorem \ref{Theorem 3} describes a situation where crack is a priori excluded, while elastic deformation is present together with possible damage due to plastic yielding on the one-dimensional subset $K$: the free gradient-discontinuity set. 

\vskip0.2cm
In the subsequent analysis we shall use the following notation for the various
contribution to the total mechanical energy:
\begin{equation}\label{Fc}
    F_r(K,v_r)\ :=\ \mathcal{H}^1\left( K \cap \overline \Omega\right)\, +\,
                \int_{\Omega\setminus K}\!  \Big( \,\eta\,|D^2 v_r|^2  -f_rv_r  \,\Big)\,d\xx \ ,
\end{equation}
\begin{equation}\label{MM}
         M(v_r-v_p)\ :=\           \mu\!\int_{\Omega}\! |v_r-v_p|^2 \, d\xx \ ,
\end{equation}
\begin{equation}\label{Fp}
        F_p(v_p)\ :=\           \!\int_{\Omega} \!\!\left(\gamma\,|D^2 v_p|^2
               -  f_p \, v_p \right) \, d\xx\ .
\end{equation}
Hence
\begin{eqnarray}\label{F=Fc+G+Fp}
\label{EE}
    & \,E(K,v)=F_r(K,v)+M(v-w) \,, \qquad\quad\ f_r=f, & \  \, \hbox{with domain }\eqref{admissible K v}, 
    \vspace{0.2cm}    
    \\ 
    \label{FF}
    & \ F(K,v_r,v_p)= F_r(K,v_r) + M(v_r-v_p) + F_p(v_p) \,,
    & \ \ \hbox{with domain }\eqref{admissible K vr vp}
     \vspace{0.2cm} 
     \\
     \label{MMM}
    & \  \ G(K,v_r,v_p)= F(K,v_r,v_p) + \sigma \int_{K\cap \overline\Omega}|\,[Dv_r]\,|\,d\mathcal H^1 \, , & \ \ \hbox{with domain }\eqref{admissible K vr vp G}
\end{eqnarray}
where
$F_r$ represents the potential energy of the reinforcement under the
Griffith assumption on the fracture energy, $M$ represents the adhesive interaction energy
(dependent on the slip $|v_r-v_p|$ between the plate and the reinforcement) and $F_p$ represents the elastic energy of the Kirchhoff-Love plate under the action of a transverse dead load $f$.
\vskip0.2cm

\begin{rmk}
We emphasize that, when minimizing \eqref{F}, the Dirichlet datum turns out to be forced on the plate ($v_p=w $ on $\partial \Omega$) since $v_p-w\!\in\! H^2_0(\Omega)$ and $H^2(\Omega_p)\!\subset \!C^0(\Omega_p),$
while the Dirichlet datum is prescribed by penalization on the reinforcement (through $v_r=w$ a.e. $\Omega_p\setminus
\overline\Omega$). Hence the damage of the reinforcement may develop also at the boundary: if this is the case then $ \mathcal{H}^1(K\cap\partial\Omega)>0$.
\\ 
In any case: $K\subset \overline\Omega$; $K$ is the closure of the set where either $v_r$ or $\nabla v_r$ is not continuous;
$w\in C^2$ and $v_r=w$ in $\Omega_p\setminus \overline\Omega$.
\end{rmk}

\begin{rmk}
The notions of essential pair or triplet in \eqref{admissible K v},\eqref{admissible K vr vp},\eqref{admissible K vr vp G} select those pairs or triplets which are cleansed of every
artifact that does not affect the functional value and are good representatives in equivalence
classes of admissible displacements. These classes allow highly irregular displacement function $v$ for $v_r$ for the reinforcement:
see Remarks 2.3-2.5 and Lemmas 2.6, 2.7 in \cite{CLTadvcv} 
for comparison with Definition 2.1
in \cite{JMPA} of
admissible triplets in the context of image segmentation and/or image inpainting.
Minimization among admissible triplets (as defined in \cite{SBZDIR}) would be equivalent
to minimization among essential admissible triplets.
\end{rmk}
\begin{rmk}\label{rmk 2.7}
The more general case where $|D^2v_p|^2$ is replaced by $Q(D^2v_p)$,
with $Q$ positive definite quadratic form, leads to claims similar to the ones we prove here (in Theorems \ref{Theorem 2} and \ref{Theorem 3}) without any change in the proofs.\\
\end{rmk}
\begin{rmk} The present paper deals with the Dirichlet boundary condition for both reinforcement and plate:
explicitly the reinforcement acts on the whole plate $\Omega_p$ and sticks perfectly to it outside $\overline{\Omega}.$\\
Nevertheless the study of Neumann boundary condition for the reinforcement,
still keeping the Dirichlet condition $w$ on the plate
(this boundary conditions correspond to a structure where the reinforcement is present only on
the proper subset $\Omega$ of the plate $\Omega_p$),
can be easily recovered by the present analysis with minor changes:
by considering admissible displacements for the reinforcement defined only in the smaller domain reference set $\Omega $ 
and replacing \eqref{admissible K v},\eqref{admissible K vr vp},\eqref{admissible K vr vp G} respectively by
\begin{equation}\label{admissible K v Neumann}
\begin{array}{l}
    \displaystyle K \ \hbox{is the smallest closed subset of } \mathbb{R}^2   \ s.t.\, 
    \displaystyle v\in C^2(\Omega\setminus K)\,;
\end{array} 
\end{equation}
\begin{equation}\label{admissible K V Neumann}
\begin{array}{l}
     v_p\!-\!w\!\in\! H^2_0(\Omega_p)\,, \
    \displaystyle K \, \hbox{ smallest closed subset of } \mathbb{R}^2 
    \ s.t.\, 
    \displaystyle v_r\!\in\! C^2(\Omega\!\setminus\! K) \,;
\end{array}
\end{equation}
\begin{equation}\label{admissible K V Neumann plastic}
\begin{array}{l}
  v_p\!-\!w\!\in\! H^2_0(\Omega_p), \,
   \! \displaystyle K \, \hbox{smallest closed subset of } \mathbb{R}^2  s.t.\, 
    \displaystyle v_r\!\in\! C^2(\Omega\!\setminus\! K)\!\cap\! C^0(\Omega) .
\end{array}
\end{equation}
All the claims in Theorems \ref{Theorem 1} and {Theorem 2} still hold true 
under these different admissible classes
of pairs and triplets. The only change to be made in the proofs amounts to refer to \cite{SBZ}
instead of \cite{SBZDIR},
to perform the analysis of partial regularity for weak minimizers.
\end{rmk}
\section{One-dimensional analysis: reinforcement of flexural beams}\label{1D}
We study the 1D case, namely the hard-device reinforcement and
strengthening reinforcement of a clamped beam, in order to make explicit
some properties of minimizers
like compliance identity, Euler equations, issues related to uniqueness and possible addition of the unilateral constraint describing the non-interpenetration of beam and reinforcement. 
The displacement of the clamped beam is modeled by a function of one variable
which is free in the interval $[-1,+1]$ while it must coincide with a given function
outside $[-1,+1]$ to take into account of boundary conditions.
%
We consider possibly different weights for energy dissipation when
crack or crease do appear: the constants $\alpha$ and $\beta$ introduced below.
\\
%
We consider real valued functions defined on bounded intervals, and set
\begin{equation}\label{ass 1-d}
  \Omega \!=\!(-1,1),\quad \Omega_p \!=\! (-2,2), \quad
  w\in C^2(-2,+2),
\end{equation}
moreover, concerning the notation, $\dot v$ denotes the absolutely continuous part of the distributional derivative $v'$ of $v$, 
$\ddot v$ denotes the absolutely continuous part of $(\dot v)'\,,$
$S_v$ denotes the set of discontinuity points of $v$, $S_{\dot v}$ denotes
the set of discontinuity points of $\dot v$ and
$v^-$,\,$v^+$ denote respectively the left and right limit of $v$.
Since we will consider admissible only piece-wise $H^2$ functions $v:(-2,2) \to \mathbb{R}$
fulfilling $v=w$ in $(-2,-1)\cup (1,2)$, we have $\left(S_v\cup S_{\dot v}\right)\subset[-1,1]$
for them all.
Here $H^2(a,b)$ denotes the usual Sobolev space of real-valued functions $v\in L^2(a,b)$ s.t. $v', v''\in L^2(a,b)$.
\\
We emphasize that the beam may develop singularities also
at both clamped endpoints $\pm 1$:
namely, it may undergo crack discontinuity 
(if $S_v\cap \{\pm 1\}$ is nonempty) or plastic-yield bending (if $S_{\dot v}\cap \{\pm 1\}$ is nonempty).
\\
After labeling by $\sharp$ the counting measure, we denote by
\begin{equation}\label{damage1D}
    J(v)\,=\,\alpha\sharp\,( S_v) +
    \beta \sharp\,(S_{\dot v}\setminus S_{ v})
\end{equation}
the whole energy associated to damage of the reinforcement: in this one-dimensional setting we allow different release energy for crack and crease,
respectively $\alpha$ and $\beta$. 
In the one-dimensional setting
the functionals $E,\,F$ and $G$ are replaced respectively by $E_1,\,F_1$ and $G_1$ defined below:
we emphasize that for them all the strong and weak formulation of  related free discontinuity problems coincide in the one-dimensional case, since finite energy entails that only a finite number of discontinuity points is allowed by finite energy, hence only piece-wise regular functions have finite energy.\\
The total energy for hard-device reinforcement of a clamped beam is given by functional $E_1$:
\begin{equation}\label{E1d}
\begin{array}{l}\displaystyle
    E_1(v)\ =\ J(v)
    \, +\,
               \int_{-1}^1 \big( \,\eta\,|\ddot v|^2 -f v\,\big)\, dx \,+\,
               \mu\,\int_{-1}^1  |v-w|^2 \, dx
               \ ;
    \end{array}
   \end{equation}
functional $E_1$ has to be minimized among the admissible functions $v$ such that
\begin{equation}\label{admissible 1d}
     v\in     \mathbf{H}^2(-2,2) := \{v:(-2,2)\to \mathbb{R}, \hbox{ s.t. }v\hbox{ is piece-wise }H^2\}\,,
\end{equation}
\begin{equation}\label{Dirichlet constraint n=1} v=w \qquad \hbox{a.e } (-2,-1)\cup(1,2)\,
\end{equation}
The total energy for strengthening reinforcement of a clamped beam is given by functional $F_1$:
\begin{equation}\label{F1d}
\begin{array}{l} \\
    F_1(v_r,v_p) \ =\ F_1(V)
    \ :=\   \\  \displaystyle 
    J(v_r)
     +\!
                \int_{-1}^1\! 
                \Big(\eta|\ddot v_r|^2 -f_rv_r\Big)dx +
               \mu\!\int_{-1}^1\! |v_r-v_p|^2 \, dx +
               \!\int_{-1}^1 \!\!\Big(\gamma\,|v_p''|^2
               -f_p \, v_p \Big)  dx \, .
               \end{array}
\end{equation}
Functional $F_1$ has to be minimized among the admissible pairs $V$ such that
\begin{equation}\label{admissible 1d F1}
     V=(v_r,v_p)\in
    \mathbf H^2(-2,2)\times H^2(-2,2) \hbox{ with } 
\end{equation}
\begin{equation}\label{Dirichlet constraint n=1 F1}
v_r=v_p=w \quad \hbox{in}\quad (-2,-1)\cup(1,2)
\end{equation}

The total energy for strengthening reinforcement of an elastic-plastic clamped beam is given by functional $G_1$:
\begin{equation}\label{G1d}
\begin{array}{l} \\
    G_1(v_r,v_p) \ =\ G_1(V)
    \ :=\   \\  \\ \displaystyle 
    \beta\ \sharp \big(S_{\dot v_r}\big) \,+\, \sigma\, \sum_{S_{\dot v_r}}\,|\,[\dot v_r]\,| \,+
     \\ \displaystyle
                \int_{-1}^1\! 
                \Big(\eta|\ddot v_r|^2 -f_rv_r\Big)dx +
               \mu\!\int_{-1}^1\! |v_r-v_p|^2 \, dx +
               \!\int_{-1}^1 \!\!\Big(\gamma\,|v_p''|^2
               -f_p \, v_p \Big)  dx \, .
               \end{array}
\end{equation}
Functional $G_1$ has to be minimized among the admissible pairs $V$ fulfilling
\begin{equation}\label{admissible 1d G1}
     V=(v_r,v_p)\in
    \Big( C^0([-2,2])\cap\mathbf H^2(-2,2)\Big)\times H^2(-2,2) 
    \,,
\end{equation}
\begin{equation}\label{Dirichlet constraint n=1 G1}
v_r=v_p=w \quad \hbox{in}\quad (-2,-1)\cup(1,2)\,.
\end{equation}

Concerning respectively \eqref{admissible 1d}, \eqref{admissible 1d F1}, \eqref{admissible 1d G1},
we recall that in all cases the finiteness of total energy implies 
respectively 
$\sharp \big(S_u\big)<+\infty$,
$\sharp\big(S_{u_r}\cup (S_{\dot u_r})\big)<+\infty$ and $S_{u_r}=\emptyset$ with $\sharp\big( S_{\dot u_r}\big)<+\infty$, 
hence
$u$ and $u_r$, are made by finitely many $H^2$ pieces.
\vskip0.1cm
%
\begin{thm}\label{Thm 3.1}
Assume \eqref{ass 1-d},\,\eqref{damage1D},\,\eqref{E1d},\,$\eta>0,$\,$\mu>0$, $f\in L^2(-1,1)$ and
\begin{equation}\label{alfabeta}
0<\beta\leq\alpha\leq 2\beta
,
\end{equation}
then the functional $E_1$ defined by \eqref{E1d} achieves a finite minimum over
functions $v$ fulfilling conditions \eqref{admissible 1d},\eqref{Dirichlet constraint n=1}.
\end{thm}
\textsl{Proof -}
After noticing that
\begin{equation}
\label{reduction E1 to BZ}
    \int_{-1}^1 \!\!\!\big( \,\mu\,  |v-w|^2 \, -f v\,\big)\, dx =
     \int_{-1}^1 \!\! \mu\Big(  v-\big(w+f/(2\mu)\big)\Big)^{\!2}dx
     \,-\!
     \int_{-1}^1 \!\!\Big( fw+f^2/(4\mu)\Big)dx
               \end{equation}
where the last summand on right-hand side is a constant,
we deduce that the functional $E_1$ is bounded from below since
all terms are nonnegative, except such constant. Thus
the claim follows 
by choosing $g=w+f/(2\mu)\in  L^2(-1,1)$ in the result of \cite{COSCIA}. \hfill$\square$
 constraint n=1.

\begin{thm}\label{Thm 3.2}
Assume\,\eqref{ass 1-d},\,\eqref{damage1D},\,\eqref{F1d},\,\eqref{alfabeta},\,$\eta\!>\!0,$\,$\mu\!>\!0$,\,$\gamma\!>\!0$,\,$f_r,f_p\in L^2(-1,1)$.

Then the functional $F_1$ defined by \eqref{F1d} achieves a finite minimum over
functions $v$ fulfilling conditions \eqref{admissible 1d F1},\eqref{Dirichlet constraint n=1 F1}.
\end{thm}
\textsl{Proof -}

The only novelty with respect to Theorem \ref{Thm 3.1} consists in the addition of
the functional 
$\int_{-1}^1 \!\!\left(\gamma\,|v_p''|^2 -  f_p \, v_p -  f_r \, v_r \right) dx $ 
and adhesive interaction $\mu\int_{-1}^1 |v_r-v_p|^2dx$
coupling $v_r$ and $v_p$\,.\\
In case of functional $F_1$ the identity \eqref{reduction E1 to BZ} reads as follows
 \begin{equation}
\label{reduction F1 to BZ}
     \int_{-1}^1 \!\!\big( \,\mu\,  |v_r-v_p|^2 \, -f_r v_r\,\big)\, dx =
     \int_{-1}^1 \!\! \mu\Big(  v_r-\big(v_p+f_r/(2\mu)\big)\Big)^2dx\,-
     \int_{-1}^1 \!\!\Big( f_r v_p+f_r^2/(4\mu)\Big)dx
               \end{equation}         
where the last summand is not a priori bounded from below,
unless we show an a priori bound on $\|v_p\|_{L^2(-1,1)}$,
moreover we have to check that minimizing sequences            
are not made by pair sequences $\big((v_r)_n,(v_p)_n\big)$ 
balancing $ \mu\|(v_p)_n\|_{L^2(-1,1)}^2\to +\infty$ together with 
$  \int_{-1}^1 \left(\gamma\,|v_p''|^2
               -  f \, v_p \right)  dx \to-\infty$. 

 This is prevented by the subsequent estimate from below 
 \eqref{F1 inf boundedness} due to $\,v_p-w\in H^2_0(-1,1)$,
where  $C_P$ denotes the best Poincar\'e constant fulfilling  
\begin{equation}\label{Poincare1}
\|v\|_{L^2(-1,1)}^2 \,\le \, C_P\, \|v''\|_{L^2(-1,1)}^2 
 \qquad \forall\,
v\!\in\! H^2_0(-1,1).
\end{equation}
%
and we denote shortly $\|\cdot\|^2$ in place of $\|\cdot\|^2_{L^2(-1,1)}$\,:
\vskip-0.3cm
\begin{equation}
\label{F1 inf boundedness}
 \begin{array}{l} \\
    F_1(v_r,v_p) \ =
      \\  \displaystyle 
    = J(v_r)
     +\!
                \int_{-1}^1\! \!
                \Big(\eta|\ddot v_r|^2 +
               \mu  |v_r-v_p|^2  
              (\gamma\,|v_p''|^2 \Big) \, dx \,
             - \! \int_{-1}^1\! \left( f_rv_r +f_p \, v_p \right)  dx  
               \\
               \displaystyle 
               = J(v_r)
     +\!
                \int_{-1}^1\! \!
                \Big(\eta|\ddot v_r|^2 +
               \mu  |v_r-v_p|^2  
              (\gamma\,|v_p''|^2 \Big) \, dx \,+
              \\
               \displaystyle 
               \qquad -  \int_{-1}^1\! \! f_r(v_r-v_p) \,dx
                -  \int_{-1}^1\! \! (f_r+f_p)(v_p-w) \,dx
                 -  \int_{-1}^1\! \! (f_r+f_p)w  \,dx \,
       \\
               \displaystyle        
        \ge\!
        J(v_r) +
                \int_{-1}^1\! \!
                \Big(\eta|\ddot v_r|^2 +
               \mu  |v_r-v_p|^2  
              +\gamma\,|v_p''|^2 \Big) \, dx \,+     
              \\
             \displaystyle   
             \ \ - \frac 1 {4\mu}\,\|f_r\|^2 \!-  \mu\,\|v_r-v_p\|^2   
             - \sqrt{C_P}\,\|f_r+f_p\|\,\|(v_p-w)''\|
             -\|f_r+f_p\| \,\|w\|     
             \\
               \displaystyle        
        \ge\!
        J(v_r) +
                \int_{-1}^1\! \!
                \Big(\eta|\ddot v_r|^2 +
               \gamma\,|v_p''|^2 \Big) \, dx \,+     
              \\
             \displaystyle   
             \quad - \frac 1 {4\mu}\,\|f_r\|^2    
             - \sqrt{C_P}\,\|f_r+f_p\|\,\|(v_p-w)''\|
             -\|f_r+f_p\| \,\|w\|   
             \\
              \displaystyle   
             \ge\!
        J(v_r) +
                \int_{-1}^1\! \!
                \Big(\eta|\ddot v_r|^2 +
               \gamma\,|v_p''|^2 \Big) \, dx \,
                     -  \sqrt{C_P}\,\|f_r+f_p\|\,\|v_p''\|
              \\
             \displaystyle   
             \quad - \frac 1 {4\mu}\,\|f_r\|^2    
             -  \sqrt{C_P}\,\|f_r+f_p\|\,\|w''\|
             -\|f_r+f_p\| \,\|w\|    
              \\
              \displaystyle   
             \ge\!
        J(v_r) +
                \int_{-1}^1\! \!
                \Big(\eta|\ddot v_r|^2 +
               \gamma\,|v_p''|^2 \Big) \, dx \,
                     - \frac {\gamma} 2 \|v_p''\|^2 +
              \\
             \displaystyle   
\ 
             - \frac {C_P} {2\,\gamma} \|f_r+f_p\|^2
             - \frac 1 {4\mu}\,\|f_r\|^2    
             -  \sqrt{C_P}\,\|f_r+f_p\|\,\|w''\|
             -\|f_r+f_p\| \,\|w\|   \,   
             \\
              \displaystyle   
             \ge\!
        J(v_r) +
                \int_{-1}^1\! \!
                \Big(\eta|\ddot v_r|^2 +
               (\gamma/2)\,|v_p''|^2 \Big) \, dx \,
                    - C(\mu,\gamma,f_r,f_p,C_P)\ .
   \end{array}
\end{equation}
 \vskip-0.3cm
Then we can fix a minimizing sequence $\big((v_r)_h,(v_p)_h\big)$ for $F_1$, and get boundedness of $\|(v_p)_h\|_{H^2(-1,1)}$,
thanks to $F_1(w,w)<+\infty$, \eqref{F1 inf boundedness} and the Poincar\'e inequality \eqref{Poincare1}. There is $u_p\in L^2$ such that we can extract a subsequence, without relabeling, fulfilling $(v_p)_h \to u_p$
weakly in $H^2$ and strongly in $L^2$, with $\|(v_p)''_h\|_{L^2}\to \|u_p''\|_{L^2}$. \\
By \eqref{reduction F1 to BZ} and \eqref{F1 inf boundedness} 
also $\|(v_r)_h-(v_p)_h\|_{L^2}$ and $\|(v_r)_h\|_{L^2}$ are bounded: by extracting again, without relabeling, $(v_r)_h \!\to\! u_r$ weakly in $L^2$.
We write
\vskip-0.6cm
\begin{equation}
\label{Fh Fh A B}
\begin{array}{l}
  F_1\big((v_r)_h,u_p\big) =         F_1\big((v_r)_h,(v_p)_h\big) +
 \Big(
  F_1\big((v_r)_h,(v_p)_k\big)  - F_1\big((v_r)_h,(v_p)_h\big)  
 \Big) +
 \vspace{-0.1cm}
 \\ 
   +
 \Big( 
   F_1\big((v_r)_h, u_p \big)   \!  - \! F_1\big((v_r)_h,(v_p)_k\big) 
 \!\Big)
 =
 F\big((v_r)_h,(v_p)_h\big) + A(h,k) + B(h,k).
 \end{array}
\end{equation}
\vskip-0.3cm
By lower semicontinuity of the functional
 $
 v \mapsto  \mu\!\int_{-1}^1\! |v-v_p|^2 \, dx +
               \!\int_{-1}^1 \!\left(\gamma\,|v''|^2
               \!-\!  f_p \, v \right)  dx
               $
 we get $\liminf_k B(h,k) \le 0$. \ Moreover, $\forall\,\varepsilon>0$ $\exists \,h_\varepsilon:$ for every $h,k>h_\varepsilon$
 \vskip-0.6cm
\begin{eqnarray*}
 |A(h,k)| 
 &\le&
 \mu\!\int_{-1}^{1}\! | (v_p)_k-(v_p)_h |  | 2(v_r)_h - ((v_p)_k+(v_p)_h) |\,dx
 +
   \\
  &&
  \ + \eta
  \big(\|(v_p)''_k\|^2-\|(v_p)''_h\|^2\big) + 
   \|f_p\|_{L^2} \big(\| (v_p)_k-(v_p)_h \| \big)
   \ < \ \varepsilon \ .
 \end{eqnarray*}
 \vskip-0.2cm
By evaluating on both sides of \eqref{Fh Fh A B} first $\liminf_k$, then $\liminf_h$, we obtain that also $\big((v_r)_h,u_p\big)$ is a minimizing sequence for the functional $F_1$ which is lower semicontinuous, or equivalently 
$(v_r)_h$ is a minimizing sequence for 
functional $E_1$
with datum  $u_p$. Then $(u_r,u_p)$ belongs to $ \argmin F_1$.             \hfill$\square$ \vskip0.2cm 

\begin{thm}\label{Thm 3.3}
Assume\,\eqref{ass 1-d},\,\eqref{damage1D},\,\eqref{G1d},\,
$\beta\!\ge\!0,\,\sigma\!>\!0,\,\eta\!>\!0,\,\mu\!>\!0,\,\gamma\!>\!0,$
and $\,f_r,\,f_p$ belong to $ L^2(-1,1).$ 
Then the functional $G_1$ defined by \eqref{G1d} achieves a finite minimum over
pairs $(v_r,v_p)$ fulfilling  
Dirichlet condition
\eqref{Dirichlet constraint n=1 F1} and
\vskip-0.5cm
\begin{equation}\label{admissible 1d G1 bis}
     V=(v_r,v_p)\in
    \Big( C^0(-2,2) \cap \mathbf H^2(-2,2)\Big)\times H^2(-2,2) 
    \,.
\end{equation}
\end{thm}
\textsl{Proof -} 
Notice that $v_r\in C^0$ entails $S_{v_r}=\emptyset.$
We set $I(v_r)\,=\,\beta\, \sharp \big(S_{\dot v_r}\big) + \sigma \sum_{S_{\dot v_r}}\,|\,[\dot v_r]\,|$.
\\
By arguing as like as in the derivation of
 \eqref{F1 inf boundedness} (the only difference consists in replacing $J(v_r)$ with $I(v_r)$), we get
 \vskip-0.5cm
\begin{equation}\label{G1 inf boundedness}
G_1(v_r,v_p) \ge
I(v_r) +
 \int_{-1}^1\! \! \Big(\eta|\ddot v_r|^2 +
               (\gamma/2)\,|v_p''|^2 \Big) \, dx \,
                    - C(\mu,\gamma,f_r,f_p,C_P)\,.
          \end{equation}
\vskip-0.3cm
If $\beta>0$, we conclude by arguing as in the proof of Theorem \ref{Thm 3.2} that a minimizing sequence has a subsequence converging to a minimum.
When $\beta=0$, after finding an optimal $u_p$ again by the argument in the proof of Theorem \ref{Thm 3.2}, we exploit Theorem 2.1 of \cite{PT1} to find the related optimal $u_r$.
\\
We emphasize that the safe load condition assumed in \cite{PT1} is unnecessary here thanks to adhesion term $\int_{-1}^1\!|v_r-v_p|^2dx$, providing boundedness from below by \eqref {G1 inf boundedness}.
\hfill $\square$
\vskip0.2cm
Next result shows that, provided the load and Dirichlet datum are suitably small, the 
strengthening reinforcement of the clamped beam ($\argmin F_1$) has a unique solution $(u_r,u_p)$ where $v_r$ has neither crack nor hinges, say $S_{u_r}\cup S_{\dot u_r}=\emptyset$.
\begin{thm}\label{Thm reg min for F1}
In addition to assumptions of Theorem \ref{Thm 3.2} we assume 
\vskip-0.6cm
\begin{equation}\label{hyp for uniq}
  \frac 1 {4\mu} 
  \|f_r\|^2_{L^2(-1,1)} + \frac {C_\Omega}{2\gamma}\|f_r+f_p\|^2_{L^2(-1,1)}
 + \gamma \|w''\|_{L^2(-1,1)}^2 + \!\int_{-1}^1\!\!\! (f_r+f_p)w\,dx
   < \beta - M  ,
\end{equation}
\vskip-0.2cm
where $\ -\infty < M=\min \big\{ F_1(v_r,v_p)\!: v_r-w\in H^2_0(-1,1),v_p-w\in H^2_0(-1,1)\big\} < +\infty. $
\\
Then there is a unique minimizer $(u_r,u_p)$ of $F_1$ over 
$\mathbf H^2\times H^2$ 
and such minimizer fulfils $S_{u_r}\cup S_{\dot u_r}=\emptyset$, \,thus \,$u_r-w\in H^2_0(-2,2)$.
\end{thm}
\textsl{Proof - } 
We denote shortly $\|\cdot\|^2$ in place of $\|\cdot\|^2_{L^2(-1,1)}$\,. First we note that
\\
$\min \big\{ F_1(v_r,v_p)\!: v_r-w\!\in\! H^2_0(-1,1),v_p-w\!\in\! H^2_0(-1,1)\big\} \le F_1(w,w) = (\eta+\gamma)\|w''\|^2 - \int_{-1}^1 (f_r+f_p)w\,dx $;

moreover, thanks to \eqref{F1 inf boundedness}, 
$F_1$ is bounded from below and coercive 
hence its infimum $M$ 
is attained and is a finite minimum.\\
Assume by contradiction that a minimizer
$(u_r,u_p)$ of $F_1$ on 
$\mathbf H^2\times H^2$ 
has $S_{u_r}\!\cup S_{\dot u_r}\!\not=\!\emptyset$, we deucee $ \alpha\,\sharp \big(S_{u_r}\big)+\beta\, \sharp \big(S_{u_r}\!\setminus\!S_{\dot u_r}\big)\ge \beta$, hence,
exploiting Poincar\'e inequality \eqref{Poincare1} and assumption \eqref{hyp for uniq} we get
 \begin{equation*}
 \begin{array}{l}
 \!\!\!\displaystyle
 \int_{-1}^1 \! f_ru_r  + \int_{-1}^1 \!f_pu_p 
 -
     \int_{-1}^1 (f_r+f_p)\,w 
=
   \int_{-1}^1 f_r\,(u_r-u_p) \,+  
   \int_{-1}^1 (f_r+f_p)\,(u_p-w)  \,\le\,
  \\ 
  \vspace{0.1cm}
 \le
 \displaystyle
   \frac 1 {4\mu} \|f_r\|^2 + \mu \|u_r-u_p\|^2 +
   \frac {C_\Omega} {2\gamma} \|f_r+f_p\|^2 + \frac {\gamma} {2C_\Omega} \|u_p-w\|^2 \le
    \\
    \vspace{0.1cm}
 \le\displaystyle
   \frac 1 {4\mu} \|f_r\|^2 + \mu \|u_r-u_p\|^2 + 
   \frac {C_\Omega} {2\gamma} \|f_r+f_p\|^2 + \frac {\gamma} 2 \, \|u_p''-w''\|^2 \le
    \\ 
 \le\displaystyle
   \frac 1 {4\mu} \|f_r\|^2  \!+ 
   \frac {C_\Omega} {2\gamma} \|f_r+f_p\|^2 
   \!+ \gamma \|w''\|^2
    +\eta \|\ddot u_r\|^2
   \!+  \mu \|u_r-u_p\|^2 \!+ \gamma \, \|u_p''\|^2
  \\
  <\displaystyle
  \beta - M 
  -
     \int_{-1}^1 (f_r+f_p)\,w 
   \,+\, \eta \|\ddot u_r\|^2
   +  \mu \|u_r-u_p\|^2 + \gamma \, \|u_p''\|^2
  \vspace{-0.2cm}
 \\
      \le\displaystyle
    \alpha\,\sharp \big(S_{u_r}\big)+\beta\, \sharp \big(S_{u_r}\!\setminus\!S_{\dot u_r}\big)
     \,+ \eta \|\ddot u_r\|^2\!
   +  \mu \|u_r-u_p\|^2\! + \gamma \, \|u_p''\|^2 \!
    - M 
    -\!
     \int_{-1}^1\!\! \!(f_r+f_p)\,w\,. 
   \end{array}
   \end{equation*}
\vskip0.1cm
say, an inequality contradicting minimality of $(u_r,u_p)$\,: \
$
 M \, <\, F_1(u_r,u_p)\,.$
\vskip0.1cm
Uniqueness of minimizer over $\mathbf H^2\times H^2$ with 
Dirichlet datum $w$
follows by uniqueness over $ H^2\times H^2$. \hfill$\square$
\begin{rmk}\label{Rmk reg min for E1}
By analogous computations to the ones in the last proof, we obtain that the inequality 
\begin{equation*}
 \label{hyp for uniq E1}
  \frac 1 {4\mu} \, \|f\|^2_{L^2(-1,1)} + \int_{-1}^1 \! f\,w\,dx 
  \ < \ \beta - \widetilde M
\end{equation*}
entails uniqueness and $H^2$ regularity
for minimizer of $E_1$ in $\mathbf H^2$
with Dirichlet boundary condition $w$, 
where $$ -\infty \,<\, \widetilde M\,=\,\min \big\{ E_1(v)\!: v-w\in H^2_0(-1,1)\big\}\, \le\, E_1(w) < +\infty\,.$$
\end{rmk}
\vskip0.2cm
We show the analysis of $E_1$ under the addition of
the unilateral constraint 
\begin{equation}\label{eq unilateral constraint}
\qquad v\ge w \qquad \hbox{ on } [-1,+1] \,.
\end{equation}
Concerning notation, from now on we set $v^+(x)=\lim_{t\to x_+} v(t)$, $v^-(x)=\lim_{t\to x_-} v(t)$.
\vskip0.1cm
\begin{rmk}\label{rmk constraint}
Actually, the constraint \eqref{eq unilateral constraint} has to be understood as a pointwise everywhere weak inequality, since it refers to functions $v\in\mathbf H^2(-2,2)$: explicitly,
$v(x)\ge w(x)$ at $x\in [-1,+1]\!\setminus\!S_v$;
$v^+(x)\ge w(x)$ at $x\in  S_v\cup\{-1\}$;
$v^-(x)\ge w(x)$ at $x\in S_v\cup\{+1\}$.
\\
Thus, the contact set $\{x\!\in\! [-2,2]\!: v^+(x)\!=\!w(x) \hbox{ or } v^-(x)\!=\!w(x)\}$ is a closed set for every $v$ fulfilling \eqref{eq unilateral constraint}; the complement in $(-1,1)$ of the contact set is an open set.
\\
Actually, the inequality \eqref{eq unilateral constraint} prevents interpenetration and refers to a reinforcement placed above: this conventional choice is made here in order to have agreement with the usual formulation of variational inequalities (\cite{BC}).
\end{rmk}
%
\begin{thm}\label{Thm 3.5}(hard device with unilateral constraint)\\
Assume \eqref{ass 1-d},\eqref{damage1D},\eqref{E1d},\eqref{alfabeta} $\eta>0,$ $\mu>0$ and $f\!\in\! L^2(-1,1)$.\\
Then the functional $E_1$ achieves a finite minimum over
pairs $(v_r,v_p)$ fulfilling conditions \eqref{admissible 1d},\eqref{Dirichlet constraint n=1} together with the unilateral constraint 
\eqref{eq unilateral constraint}.
\end{thm}
\textsl{Proof - } The proof can be achieved by exact repetition of the 
argument in the proof of Theorem \ref{Thm 3.1} for the unconstrained case: both unilateral constraint $v\ge w$ on $[-1,+1]$ and Dirichlet condition $v= w$ a.e. on $(-2,-1)\cup(1,2)$ affect neither the compactness, nor the lower semicontinuity properties of $E_1$; moreover the a.e. convergence  preserves the constraint in the limit of minimizing sequences. \hfill $\square$
\vskip0.1cm
\begin{thm}\label{Thm reinf1 with unilat constr}(reinforcement with unilateral constraint)\\
Assume
\,\eqref{ass 1-d},\,\eqref{damage1D},\,\eqref{F1d},\,\eqref{alfabeta},\,$\eta\!>\!0,$\,$\mu\!>\!0$,\,$\gamma\!>\!0$,\,$f_r,f_p\in L^2(-1,1)$.
\\
Then the functional $F_1$ achieves a finite minimum over
pairs $(v_r,v_p)$ fulfilling the conditions 
\eqref{admissible 1d F1} and \eqref{Dirichlet constraint n=1 F1} together with the unilateral constraint (corresponding to a reinforcement placed above the plate)
\begin{equation}\label{eq unilateral constraint F1}
\qquad v_r\ge v_p \qquad \hbox{ on } [-1,+1] \,.
\end{equation}
\end{thm}
Also the constraint \eqref{eq unilateral constraint F1} has to be understood as a pointwise  everywhere everywhere weak inequality,
in the sense of Remark \ref{rmk constraint}, as like as \eqref{eq unilateral constraint} but here with $v_p$ replacing $w$: thus, the admissible pairs belong to the convex set 
\begin{equation*}
\mathbf K \!:=\! \left\{(v_r,v_p)\!\in\! \mathbf H^2(-2,2)\!\times\! H^2(-2,2)\!:\, v_r\!\ge\! v_p \,\hbox{on} \, [-1,1], \,  v_r \!=\!v_p\!=\! w \, \,\hbox{on}\, (-2,-1)\!\cup\!(1,2)\right\}
\end{equation*}
\vskip0.1cm
\textsl{Proof of Thm \ref{Thm reinf1 with unilat constr} -}
The proof can be achieved by exact repetition of the 
argument in the proof of Theorem \ref{Thm 3.2} for the unconstrained case: both unilateral constraint $v_r\ge v_p$ a.e. on $(-1,+1)$ and Dirichlet condition $v_r\!=\!v_p\!=\! w$ a.e. on $(-2,-1)\cup(1,2)$ affect neither the compactness, nor the lower semicontinuity properties of $F_1$; moreover a.e. convergence  preserves the constraint in the limit of minimizing sequences. \hfill $\square$
\vskip0.1cm
By performing all the admissible variations of minimizers for $E_1$ and $F_1$ without the unilateral constraint, we can deduce the necessary conditions for minimality listed below in Propositions \ref{EulerThm} and \ref{ComplianceThm}.
\begin{prop}\phantom{.}\label{EulerThm}\vspace{-0.2cm}
(\textbf{Euler equations for functional }$E_1$) Every $u\in\argmin E_1$ fulfils
\begin{equation}\label{EulerE1d}
    \eta \,u''''+\mu (u-w)=f/2 \qquad \hbox{in } (-1,1)\setminus
     \left( S_{u}\cup S_{\dot u} \right) \ ,
\end{equation}
\begin{equation}\label{EulerE1dSu}
    \ddot u^+ = \ddot u^- = \dddot u^+= \dddot u^-  = 0 \qquad\qquad\hbox{in }
     S_{u}\setminus\{\pm 1\}\ ,
\end{equation}
\begin{equation}\label{EulerE1dSDu}
    \ddot u^+ = \ddot u^- = [\dddot u]  = 0 \qquad\qquad\ \hbox{on }
     S_{\dot u}\setminus\left( S_u\cup \{\pm 1\} \right) \ ,
\end{equation}
\begin{equation}\label{EulerE1dsintesi} \ddot u\in H^2(-1,1)\qquad\hbox{and}\quad \eta\,(\ddot u)''+\mu (u-w)=f/2 \quad \hbox{on }\mathcal{D}'(-1,1)\ .
     \end{equation}
\begin{equation}\label{EulerE1endpointsSu}
\left\{\
\begin{array}{l}
    \ddot u^+ (-1)= \dddot u^+(-1) =  0 \qquad \hbox{if }-1
    \in S_u\setminus S_{\dot u}\ , \vspace{0.2cm} \\
    \ddot u^- (+1)= \dddot u^-(+1) =  0 \qquad \hbox{if }+1 \in S_u\setminus S_{\dot u}\  ,
    \end{array} \right.
\end{equation}
\begin{equation}\label{EulerE1endpointsSDumenoSu}
\left\{\
\begin{array}{l}
    \ddot u^+ (-1)=   0 \qquad \hbox{if }-1\in S_{\dot u}\setminus S_u\ , \vspace{0.2cm} \\
    \ddot u^- (+1)=   0 \qquad \hbox{if }+1\in S_{\dot u}\setminus S_u\  ,
    \end{array} \right.\qquad\phantom{.}
   \end{equation}
When $\alpha\!=\!\beta$ the conditions \eqref{EulerE1dSu},\eqref{EulerE1dSDu},\eqref{EulerE1endpointsSu},\eqref{EulerE1endpointsSDumenoSu} altogether are improved as follows:
\begin{equation}\label{EulerE1dSuimproved}
  \hbox{if }\alpha=\beta \quad\hbox{ then } \quad
  \ddot u^+ = \ddot u^- = \dddot u^+ = \dddot u^- = 0 \quad
    \hbox{ on }\left( S_{u}\cup S_{\dot u} \right)  .
\end{equation}
\vskip0.1cm
(\textbf{Euler equations for functional }$F_1$)
Every $(u_r,u_p)\in\argmin F_1$ fulfils
\begin{equation}\label{EulerF1dprima}
    \eta\,u_r''''+\mu (u_r-u_p)=f_r/2 \qquad (-1,1)\!\setminus\!
    \left( S_{u_r}\cup S_{\dot u_r}\right)\ ,
\end{equation}
\begin{equation}\label{EulerF1dseconda}
    \gamma\, u_p''''+\mu (u_p-u_r)=f_p/2 \qquad \hbox{ in }\mathcal{D}'(-1,1)\ ,
\end{equation}
\begin{equation}\label{EulerF1dterza}
    \eta\,u_r''''+\gamma u_p''''=(f_r+f_p)/2 \qquad (-1,1)\!\setminus\!
     \left(S_{u_r}\cup S_{\dot u_r}\right)\ ,
\end{equation}
\begin{equation}\label{EulerF1dSu}
    \ddot u_r^+ = \ddot u_r^- = \dddot u_r^+ = \dddot u_r^- = 0 \qquad\qquad\
     S_{u_r}\setminus\{\pm1\}\ ,
\end{equation}
\begin{equation}\label{EulerF1dSDu}
    \ddot u_r^+ = \ddot u_r^- = [\dddot u_r]  = 0 \quad\qquad
     S_{\dot u_r}\!\setminus\!\left( S_{u_r}\cup \{\pm 1\} \right) \ ,
\end{equation}
hence $\ddot u_r\in H^2(-1,1)$ and 
\begin{equation}\label{EulerF1dsintesi1} \eta\,(\ddot u_r)''+\mu (u_r-u_p)=f_r/2 \qquad \hbox{ in }\mathcal{D}'(-1,1)\ ,
     \end{equation}
     \begin{equation}\label{EulerF1dsintesi2} \eta\,(\ddot u_r)''+\gamma u_p''''=(f_r+f_p)/2\qquad \hbox{ in }\mathcal{D}'(-1,1)\ .
     \end{equation}
\begin{equation}\label{EulerF1endpointsSu}
\left\{\
\begin{array}{l}
    \ddot u_r^+ (-1)= \dddot u_r^+(-1) =  0 \qquad \hbox{if }-1
    \in S_{u_r}\!\setminus\! S_{\dot u_r}\ , \vspace{0.2cm} \\
    \ddot u_r^- (+1)= \dddot u_r^-(+1) =  0 \qquad \hbox{if }+1
   \in S_{u_r}\!\setminus\! S_{\dot u_r}\  ,
    \end{array} \right.
\end{equation}
\begin{equation}\label{EulerF1endpointsSDumenoSu}
\left\{\
\begin{array}{l}
    \ddot u_r^+ (-1)=   0 \qquad \hbox{if }-1\in S_{\dot u_r}\!\setminus\! S_{u_r}\ , \vspace{0.2cm} \\
    \ddot u_r^- (+1)=   0 \qquad \hbox{if }+1\in S_{\dot u_r}\!\setminus\! S_{u_r}\  ,
    \end{array} \right.\qquad\phantom{.}
   \end{equation}
     So, \eqref{EulerF1dsintesi1} and \eqref{EulerF1dsintesi2}, together give
\begin{equation}\label{EulerF1dsintesi3} 2 \,\eta\,(\ddot u_r)''+\mu (u_r-u_p)+\gamma u_p''''= f_r + (1/2)f_p
\qquad \hbox{ in }\mathcal{D}'(-1,1)\ .
     \end{equation}
When $\alpha\!=\!\beta$ the conditions 
\eqref{EulerF1dSu},\eqref{EulerF1dSDu},\eqref{EulerF1endpointsSu},\eqref{EulerF1endpointsSDumenoSu} altogether are improved as follows:
\begin{equation}\label{EulerF1dSuimproved}
  \hbox{if }\alpha=\beta \quad\hbox{ then } \quad
  \ddot u_r^+ = \ddot u_r^- = \dddot u_r^+ = \dddot u_r^- = 0 \quad
    \hbox{ on }\left( S_{u_r}\cup S_{\dot u_r} \right)  .
\end{equation}
\end{prop}
Eventually we deduce the following compliance identities.
\begin{prop}\vskip0.1cm\noindent\label{ComplianceThm}
\textbf{Compliance identity for functional $E_1$:}\\ Assume $w(-1)=w(1)=w'(-1)=w'(1)=0$. Then any $u\in\argmin E_1$ fulfils  \begin{equation}\label{ComplianceE1d} E_1(u) \,=\,
J(u)
\,+\,\mu\int_{-1}^1 (\,w^2\,-\,wu\,)\,dx \
-\ \frac 1 2 \int_{-1}^1 
f u
\,dx\ .
     \end{equation}
If boundary conditions are nonhomogeneous
then the right-and side of compliance \eqref{ComplianceE1d} has to be added with the correction
\,$+\,\eta\, \left[\, \ddot u \,\dot u - \dddot u \,u\,\right]_{-1}^{+1}$\,, where
\,$[z]_{-1}^1=z^-(1)-z^+(-1)$.
Notice that (due to \eqref{EulerE1endpointsSu},\eqref{EulerE1endpointsSDumenoSu},\eqref{EulerE1dSuimproved}) some of the four terms in the correction may be null if one endpoint or the other belongs to $S_{\dot u}\cup S_u$.

     \vskip0.2cm
\textbf{Compliance identity for functional $F_1$:}
\\ Assume $w(-1)=w(1)=w'(-1)=w'(1)=0$, then any $(u_r,u_p)\in\argmin F_1$ fulfils
\begin{equation}\label{ComplianceF1d} F_1(u_r,u_p) \,=\,
J(u_r)
\,-\, \mu \int_{-1}^1  (u_r-u_p)^2  dx
\,\,-\,\frac 1 2 \,\int_{-1}^1 
f_r\, u_r
\,dx\,
\,-\,\frac 1 2 \,\int_{-1}^1 
f_p\, u_p
\,dx\ .
     \end{equation}
If boundary conditions are nonhomogeneous
then the right-and side of compliance \eqref{ComplianceF1d} has to be added with the correction

$+\left[\,\gamma\,(u_p''w'\!-u_p'''w )+ \eta\,(\,\ddot u_r\, \dot u- \dddot u_r\, u)\,\right]_{-1}^{+1}$\,, where
\,$[z]_{-1}^1=z^-(1)-z^+(-1)$.

     \end{prop}

\vskip0.3cm
\textsl{Proof of Proposition \ref{EulerThm} (Euler equations for $E_1$ and $F_1$)} -\\
Let $u$ be a minimizer of $E_1$ among $v\in\mathbf{H}^{2}:=\mathbf{H}^{2}(-2,2)$ fulfilling
\eqref{admissible 1d} and \eqref{Dirichlet constraint n=1}. For
any $v\in \mathbf{H}^2$ we set $\mathbf{\big[\kern-5pt\big[}v\mathbf{\big]\kern-5pt%
\big]}=v^{+}-v^{-}$ where $v^{-},\ v^{+}$ denote respectively the left and
right values of $v$ on $S_{v}$.\newline
We introduce the localized version of functional $E_1$:
given $w$, $\alpha $, $\beta $, we set, for any $v$ in $\mathbf{H}^2$
and any Borel set $A\subset \lbrack -1,1]$,
\begin{multline}
E_1(v,A)\ = \
\int_{A}\big(\eta|\ddot{v}^{2}|+\mu|v-w|^{2}\big)\,dx+\alpha \,\sharp(S_{v}\cap A)+\beta \,\sharp\left( (S_{\dot{v}}\setminus
S_{v})\cap A\right)\,.   \label{formula 1}
\end{multline}

\textbf{Step 1 - (Green formula)} Assume:
$u\in\argmin E_1$.\\
Since $J(u)\leq E_1(u)<+\infty$,
the set $S_{u}\cup S_{\dot{u}}$ is finite and contained in $[-1,1]$;
$u\in H^4(I)$ for every interval $I\subset (-2,2)\setminus \{S_{u}\cup S_{\dot{u}}\}$.\\
From now on we label
$t_0=-1$ and $t_{T+1}=1$ and $t_{j}$, for
$j=1,...,\mathsf{T}$, the (possibly empty) finite ordered set $\left(S_{u}\cup S_{\dot{u}}\right)\cap (-1,+1)$. 
Then, integrating by parts, the
next identity is achieved for every $\varphi \in \mathbf{H}^{2}$ 
\begin{multline}
 \label{formula 2}
\sum_{l=0}^{\mathsf{T}}\int_{t_{l}}^{t_{l+1}}\!\ddot u\,\ddot\varphi
\ dx
\ =\
\sum_{l=0}^{\mathsf{T}}\int_{t_{l}}^{t_{l+1}}\!u^{\prime \prime }\varphi
^{\prime \prime }\,dx
\ =\ 
\sum_{l=0}^{\mathsf{T}}\int_{t_{l}}^{t_{l+1}}u^{\prime
\prime \prime \prime }\varphi \,dx\ +
\\
\sum_{l=1}^{\mathsf{T}}\Big(\big(
-\dddot u ^-
(t_{l+1})
\varphi ^{-}(t_{l+1})+
\dddot u^+
(t_{l})
\varphi^{+}(t_{l})\big)+
\big(
\ddot u^-
(t_{l+1})
\dot \varphi^-
(t_{l+1})-
\ddot u^+
(t_{l})
\dot \varphi^+
(t_{l})
\big)\Big)\,+
\\
+\Big(
- \dddot u^-(t_1)\varphi^-(t_1)\,+\,\ddot u^-(t_1)\dot\varphi^-(t_1)
\,+\,
\dddot u^+(t_T)\varphi^+(t_T)\,-\,\ddot u^+(t_T)\dot\varphi^+(t_T)
\Big)\,
+
\\
+\Big(
\dddot u^+(-1)\varphi^+(-1)\,-\,\ddot u^+(-1)\dot\varphi^+(-1)
\,-\,
\dddot u^-(1)\varphi^-(1)\,+\,\ddot u^-(1)\dot\varphi^-(1)
\Big)\,.
\end{multline}
\textbf{Step 2 -} At first we show that each minimizer $u$ solves the fourth
order elliptic equation \eqref{EulerE1d} on the interior of $(-1,1)\setminus (S_{u}\cup
S_{\dot{u}})$, by performing smooth variations. For every open set $A\subset
\subset (-1,1)\setminus (S_{u}\cup S_{\dot{u}})$, for every $\varepsilon \in
\mathbb{R}$ and for every $\varphi \in C_{0}^{\infty }(A)$ we have
\begin{equation*}
0\leq E_1(u+\varepsilon \varphi ,A)-E_1(u,A)=2\varepsilon \bigg(%
\eta\!\int_{A}u^{\prime \prime }\varphi ^{\prime \prime }\,dx+\mu\!\int_{A}(u-w)\varphi
\,dx\,-\!\int_A \frac f 2 \varphi\,dx \bigg)+o(\varepsilon )
\end{equation*}
where $o(\varepsilon )$ is an infinitesimal of higher order than $%
\varepsilon $. Hence
\begin{equation*}
\eta\,\int_{A}u^{\prime \prime }\varphi ^{\prime \prime
}\,dx\ =\int_{A} \Big(f/2-\,\mu\,(u-w)\Big)\varphi \,dx\qquad
\forall\varphi \in C_{0}^{\infty }(A)
\ .
\end{equation*}
Then \eqref{EulerE1d} follows integrating
by parts with Green formula \eqref{formula 2}.\\
Now we seek the Euler conditions at inner discontinuity points and at clamped endpoints.
\vskip0.1cm
\textbf{Step 3 -} We prove necessary conditions \eqref{EulerE1dSu}
for extremality on $S_{u}\,$ and necessary conditions \eqref{EulerE1endpointsSu}
for extremality at endpoints when they do belong to $S_u\!\setminus\!S_{\dot u}$.\\

Choose $\varphi \in \mathbf{H}^{2}\cap C^{2}([t_{l},t_{l+1}])$, $%
l=0,...,\mathsf{T}$, $\mbox{\textnormal{spt}}(\varphi )\subset A$, where $A$
is a Borel subset of $[-1,1]$ with $(S_{\dot{u}}\setminus S_{u})\cap A=\emptyset $. Then
for every $\varepsilon \in \mathbb{R}$ we have
\vskip-0.5cm
\begin{equation*}
\big(S_{u+\varepsilon \varphi }\cup S_{\dot{u}+\varepsilon \dot{\varphi}}%
\big)\cap A\ \subset \ {S_{u}}\cap A
\end{equation*}
By (\ref{formula 2}) we have:
\vskip-0.5cm
\begin{eqnarray*}
\vspace{-0.2cm}
0  & \leq &
   E_1(u+\varepsilon \varphi ,A)-E_1(u,A) \ =
\\
   & = & \alpha \left( \sharp (S_{u+\varepsilon \varphi }\cap A)-\sharp (S_{u}\cap
A)\right) +\beta \sharp \left( (S_{\dot{\varphi}}\setminus S_{u+\varepsilon
\varphi })\cap A\right) \,+
\\
   &  & \qquad + \
   2\varepsilon \left( \sum_{l=0}^{\mathsf{T}}\int_{t_{l}}^{t_{l+1}}\big(\eta\, u^{\prime
\prime }\varphi ^{\prime \prime }+\mu(u-w)\varphi - \frac f 2 \varphi\big)\,dx\right) \, + \, o(\varepsilon )\, =
\\
   & = &
   \alpha \left( \sharp (S_{u+\varepsilon \varphi }\cap A)-\sharp (S_{u}\cap
A)\right) +\beta \sharp \left( (S_{\dot{\varphi}}\setminus S_{u+\varepsilon
\varphi })\cap A\right) \,+
\\
   &  & \qquad + \
   2\varepsilon \bigg(\sum_{l=0}^{\mathsf{T}}\int_{t_{l}}^{t_{l+1}}\big(\eta\, u^{\prime
\prime \prime \prime }\varphi +\mu(u-w)\varphi - \frac f 2 \varphi\big)\,dx \,+
\\
   &  & \qquad + \ {\dddot{u}}^{+}(-1)\varphi ^{+}(-1)-\ddot{u}^{+}(-1)\dot{\varphi}^{+}(-1)-{%
\dddot{u}}^{-}(1)\varphi ^{-}(1)+\ddot{u}^{-}(1)\dot{\varphi}^{-}(1) \,+
\\
   &  & \qquad + \ \,\eta\!\sum_{\left(S_{u}\cap A\right)\setminus\{\pm 1\}}\!\Big(\mathbf{\big[\kern-5pt\big[}\dddot{u}\varphi \mathbf{%
\big]\kern-5pt\big]}-\mathbf{\big[\kern-5pt\big[}\ddot{u}\dot{\varphi}%
\mathbf{\big]\kern-5pt\big]}\Big)\bigg)\ +\ o(\varepsilon )\,.
\end{eqnarray*}
Up to a finite set of possible values of $\varepsilon $ entailing cancellation of discontinuity, we have $S_{u+\varepsilon
\varphi }\cap A=S_{u}\cap A$. Then by discarding such values we can choose arbitrarily small $%
\varepsilon $ satisfying
\vskip-0.6cm
\begin{equation*}
\vspace{-0.1cm}
\sharp \,((S_{\dot{\varphi}}\setminus S_{u+\varepsilon \varphi })\cap
A)=\sharp \,((S_{\dot{\varphi}}\setminus S_{u})\cap A)=0
\end{equation*}
By taking into account \eqref{EulerE1d} and the arbitrariness of the two traces of $\varphi
$ and $\dot{\varphi}$ on the two sides of points in $S_{u}$, for small $%
\varepsilon $, we can choose $\varphi $ with $\varphi ^{\pm }=0$, and $\dot{%
\varphi}^{+}=0$ together with $\dot{\varphi}^{-}$ arbitrary, or viceversa to get $\,\ddot u^\pm=0$ on $S_u\setminus\{\pm 1\}$.\\ Similarly, we obtain $\,\dddot u^\pm=0$ on $S_u\setminus\{\pm 1\}$ by choosing $\dot{%
\varphi}^{\pm }=0$, and $\varphi ^{+}=0$ together with $\varphi ^{-}$ arbitrary or
vice-versa. So \eqref{EulerE1dSu} is proved.\\
If some clamped endpoint ($-1$ and/or $+1$) belong to $S_u$, then \eqref{EulerE1endpointsSu} is obtained as above, but taking into account that $\varphi\equiv 0$ outside $[-1,1]$.
\vskip0.1cm
\textbf{Step 4 -} We prove the necessary condition \eqref{EulerE1dSDu} for extremality on
$S_{\dot{u}}\setminus \left(S_u\cup \{\pm 1\}\right)$:
\begin{equation}
\begin{array}{cc}
\ddot{u}^{\pm }=0 & \text{in }S_{\dot{u}}\setminus\left(S_u\cup \{\pm 1\}\right)\ ,
\end{array}
\label{formula 5}
\end{equation}
\begin{equation}
\begin{array}{cc}
\mathbf{\big[\kern-5pt\big[}\dddot{u}\mathbf{\big]\kern-5pt\big]}=0 & \text{%
in }S_{\dot{u}}\setminus\left(S_u\cup \{\pm 1\}\right) \ .
\end{array}
\label{formula 6}
\end{equation}
Let $\varphi \in \mathbf{H}^{2}\cap C^{2}([t_{l},t_{l+1}])$, $l=0,...,%
\mathsf{T}$, $\mbox{\textnormal{spt}}(\varphi )\subset A$, with $A$ Borel subset of $(-1,1) $ and $S_{\varphi
}=\emptyset =(S_{u}\setminus S_{\dot{u}})\cap A$. Then, up to a finite set of possible values of $\varepsilon $ entailing cancelation of $\dot u$ discontinuity, we can choose $\varepsilon $ arbitrarily
small such that
\vskip-0.5cm
\begin{equation*}
(S_{u+\varepsilon \varphi }\cup S_{\dot{u}+\varepsilon \dot{\varphi}})\cap
A=S_{\dot{u}+\varepsilon \dot{\varphi}}\cap A=S_{\dot{u}} \ .
\end{equation*}
Moreover, by Green formula (\ref{formula 2}):
\begin{eqnarray*}
0 &\leq &
E_1(u+\varepsilon \varphi ,A)-E_1(u,A) \,\leq
\\
  & \leq &\beta \big(\sharp \,(S_{\dot{u}+\varepsilon \dot{\varphi}}\cap
A)-\sharp \ (S_{\dot{u}}\cap A)\big)\,+
\\
  && \qquad +\ 2\varepsilon\, \bigg(\sum_{l=0}^{\mathsf{T}}\int_{t_{l}}^{t_{l+1}}\Big(\eta\, u^{\prime
\prime }\varphi ^{\prime \prime }+\mu(u-w)\varphi - \frac f 2 \varphi\Big)\,dx\bigg)\,+\,o(\varepsilon ) \,=
\\
  & = &
  2\varepsilon \, \bigg(\sum_{l=0}^{\mathsf{T}}\int_{t_{l}}^{t_{l+1}}\Big( \eta\, u^{\prime
\prime \prime \prime }\varphi \,dx+\mu (u-w)\varphi - \frac f 2 \varphi\Big)\,dx\,+
\\
&& \qquad
+\ {\dddot{u}}^{+}(-1)\varphi ^{+}(-1)-\ddot{u}^{+}(-1)\dot{\varphi}^{+}(-1)-{%
\dddot{u}}^{-}(1)\varphi ^{-}(1)+\ddot{u}^{-}(1)\dot{\varphi}^{-}(1)\,+
\\
&& \qquad
+\ \eta \sum_{\left(S_{\dot{u}}\cap A\right)\setminus\{\pm 1\}}\Big(\mathbf{\big[\kern-5pt\big[}+\dddot{u}\varphi
\mathbf{\big]\kern-5pt\big]}-\mathbf{\big[\kern-5pt\big[}\ddot{u}\dot{\varphi%
}\mathbf{\big]\kern-5pt\big]}\Big)\bigg)\ +\ o(\varepsilon )\ .
\end{eqnarray*}
\vskip-0.2cm
By taking into account \eqref{EulerE1d}, for small $\varepsilon $ and by the
arbitrariness of $\varphi $ and of the two traces of $\dot{\varphi}$ on the
two sides of $S_{\dot{u}}$, we can choose $\varphi $ with $\varphi ^{\pm }=0$%
, and arbitrary $\dot{\varphi}^{+}=\dot{\varphi}^{-}$, to get (\ref{formula
5}).
\\ On the other hand, by choosing $\dot{\varphi}^{\pm }=0$ together with
arbitrary $\varphi$ and taking into account that $\mathbf{\big[\kern
-5pt\big[}{\varphi}\mathbf{\big]\kern-5pt\big]}=0$, we
obtain (\ref{formula 6}).\newline
Then \eqref{EulerE1dSDu} follow from (\ref{formula 5}) and (\ref{formula 6}).\\
%
\textbf{Step 5 -} The analysis of minimizers at $\left(S_{\dot u}\setminus S_u\right)\cap \{\pm 1\}$ can be done exactly in the same way as in Step 5, but taking into account that $u=w$ and $\varphi=0$ on $[-2,-1]\cup[1,2]$, thus obtaining \eqref{EulerE1endpointsSu} and \eqref{EulerE1endpointsSDumenoSu}.
\vskip0.1cm
\textbf{Step 6 -}
\eqref{EulerE1dsintesi} is a straightforward consequence of \eqref{EulerE1d}-\eqref{EulerE1dSDu}.
\vskip0.1cm
\textbf{Step 7 -} Eventually, under the additional condition $\alpha=\beta$,
we prove the refinement \eqref{EulerE1dSuimproved} of \eqref{EulerE1dSu},\eqref{EulerE1dSDu},\eqref{EulerE1endpointsSu},\eqref{EulerE1endpointsSDumenoSu} on $ \left(S_{\dot u} \cup S_u\right) $ 
for every minimizer $u$.\\
We are left only to show that
\vskip-0.6cm
\begin{equation}\label{improvement}
\ \text{ if }\alpha
=\beta \ \
\text{ then\,: } \ \ \
1\in S_{\dot u} \setminus S_u \Rightarrow \dddot u^-(1)=0\,;
\ \ 
-1\in S_{\dot u}\setminus S_u \Rightarrow \dddot u^+(-1)=0
\,.
\end{equation}
Fix a Borel set $A$ s.t. $A\subset \subset (-2,2)$, $S_{u}\cap A=\emptyset
\neq S_{\dot{u}}\cap A$.\newline
Let $\varphi \in \mathbf{H}^{2}\cap C^{2}([t_{l},t_{l+1}])$, $l=0,...,%
\mathsf{T}$ and
\vskip-0.5cm
\begin{equation*}
S_{\dot{u}}\cap A \,=\, S_{\varphi }\cap A  \qquad \hbox{ and } \qquad S_{u}\cap A \, =\,  S_{\dot{\varphi}}\cap A
\,=\, \, \emptyset \,.
\end{equation*}
Then, for every value of $\varepsilon \in\R$ we have $S_{u+\varepsilon
\varphi }\cap A=S_{\varphi }\cap A$ and
\begin{equation*}
\left( S_{u+\varepsilon \varphi }\cup S_{(\dot{u}+\varepsilon \dot{\varphi}%
)}\right) \cap A=S_{\dot{u}}\cap A\ .
\end{equation*}
By (\ref{formula 2}), \eqref{EulerE1d}, \eqref{EulerE1dSu} and \eqref{EulerE1dSDu} we have
\begin{eqnarray*}
0 \!\!&\leq&\!\! E_1(u+\varepsilon \varphi ,A)-E_1(u,A)
\\
\!\!&=&\!\!
\alpha \,\sharp \left( S_{u+\varepsilon \varphi }\cap A\right) +\beta
\left( \sharp \left( (S_{(\dot{u}+\varepsilon \dot{\varphi})}\setminus
S_{u+\varepsilon \varphi })\cap A\right) -\beta \sharp \left( S_{\dot{u}%
}\cap A\right) \right)  
\\ \!\!&&\!\!
\quad +\,2\varepsilon\, \left( \sum_{l=0}^{\mathsf{T}}\int_{t_{l}}^{t_{l}+1}\Big(\eta\, u^{\prime
\prime }\varphi ^{\prime \prime }+\mu(u-w)\varphi - \frac f 2 \varphi\Big)\,dx\right) \,+\,o(\varepsilon )
\\
&=&\!\!
\alpha \,\sharp (S_{\varphi }\cap A)+\beta \,\sharp \left( (S_{\dot{u}%
}\setminus S_{\varphi })\cap A\right) -\beta \,\sharp (S_{\dot{u}}\cap A) \\
&&\!\!
\quad +\,2\varepsilon\, \left( \sum_{l=0}^{\mathsf{T}}\int_{t_{l}}^{t_{l}+1}\Big(\eta\, u^{\prime
\prime \prime \prime }\varphi +\mu(u-w)\varphi - \frac f 2 \varphi\Big)\,dx\
+\eta\!\!\!\!\sum_{\left(S_{\dot{u}}\cap A\right)\setminus \{\pm 1\}}
\!\!\Big(%
\mathbf{\big[\kern-5pt\big[}\dddot{u}\varphi \mathbf{\big]\kern-5pt\big]}-%
\mathbf{\big[\kern-5pt\big[}\ddot{u}\dot{\varphi}\mathbf{\big]\kern-5pt\big]}%
\Big)\right.
\vspace{-0.5cm}
 \\
&&\!\!
\quad + \ \eta\,\Big({\dddot{u}}^{+}(-1)\varphi ^{+}(-1)-\ddot{u}^{+}(-1)\dot{\varphi}^{+}(-1)-{%
\dddot{u}}^{-}(1)\varphi ^{-}(1)+\ddot{u}^{-}(1)\dot{\varphi}^{-}(1)
\Big)\Bigg)
\vspace{0.2cm}
\\
&&\!\! \quad +\ o(\varepsilon )
\vspace{0.1cm}\\
&=&\!\!
\alpha \,\sharp (S_{\varphi }\cap A)-\beta \,\sharp (S_{\dot{u}}\cap
A)+2\eta\,\varepsilon \sum_{\left(S_{\dot{u}}\cap A\right)\setminus \{\pm 1\}}\mathbf{\big[\kern-5pt\big[}\dddot{u}%
\varphi \mathbf{\big]\kern-5pt\big]}+ o(\varepsilon )\ .
\end{eqnarray*}
Since $S_{\varphi }\cap A=S_{\dot{u}}\cap A$, when $\alpha >\beta $ then the
inequality is fulfilled for $\varepsilon $ small enough, hence we do not
obtain further information (recall that the necessary condition for semicontinuity $\alpha\geq \beta$ is always assumed). On the other hand, when $\alpha =\beta $, we get
\begin{multline*}
0\leq E_1(u+\varepsilon \varphi ,A)-E_1(u,A)\ = \\
=\ 2\eta\,\varepsilon \left(
\sum_{S_{\dot{u}}\cap
A}\mathbf{\big[\kern-5pt\big[}\dddot{u}\varphi \mathbf{\big]\kern-5pt\big]}
+
{\dddot{u}}^{+}(-1)\varphi ^{+}(-1)-{
\dddot{u}}^{-}(1)\varphi ^{-}(1)
\right)
+o(\varepsilon ) \,.
\end{multline*}
So the coefficient of $2\varepsilon $ must vanish, and by the
arbitrariness of the two traces of $\varphi $ at points in $S_{\dot u}\cap A$, of the right trace at $-1$
and of the left trace at $+1$, taking into account that $\varphi\equiv 0$ outside $[-1,1]$ we get \eqref{improvement}.
\vskip0.1cm
\textbf{Step 8 -}
We make explicit all the details for $E_1$ only, since
the proof of Euler equations for $F_1$ is identical. In fact $F_1(v_r,v_p)-E_1(v_r)=\int_{-1}^1(\gamma|\ddot v_p|^2-fv_p)$ is a classical integral functional: so the analysis of any minimizer $U=(u_r,u_p)$ of $F_1$
can be done by performing all the admissible variations separately for $u_r$ and $u_p$.
\hfill $\square$
 \vskip0.4cm
\textsl{Proof of Proposition \ref{ComplianceThm} (compliance identities) -}
Assume $u\!\in\!\argmin E_1$ and label
$t_0=-1$ and $t_{T+1}=1$ and $t_{j}$, for
$j=1,...,\mathsf{T}$, the (possibly empty) finite ordered set $\left(S_{u}\cup S_{\dot{u}}\right)\cup (-1,+1)$.\\ Then, by \eqref{EulerE1d}-\eqref{EulerE1endpointsSDumenoSu}, integrating by parts
on the intervals $[t_j,t_{j+1}]$ we get
\begin{eqnarray*}
\eta \int_{-1}^1|\ddot u|^2\,dx &=&
\eta \sum_{j=0}^{T}\int_{t_j}^{t_{j+1}}|\ddot u|^2\,dx  \,=\,
\\
&=&
-\,\eta \int_{-1}^1 (\ddot u)'\,\dot u\,dx \,+\,\eta\,\left[ \,\ddot u \,\dot u\,\right]_{-1}^{+1}\,=
\\
&=& \eta \int_{-1}^1 (\ddot u)''\,u\,dx \,+\,\eta\,
\left[ \,\ddot u \,\dot u\,- \,\dddot u\, u \,\right]_{-1}^{+1}\,=
\\
                                &=& \mu\,\int_{-1}^1 (w-u)\,u\,dx
                                \,+\, \int_{-1}^1 (f/2)\,u\,dx
                                \,+\,\eta 
        \left[ \,\ddot u \,\dot u\,- \,\dddot u\, u \,\right]_{-1}^{+1}\ ,
\end{eqnarray*}
here above and in the sequel, the notation $[z]_{-1}^{+1}$ stands for $z^-(1)-z^+(-1)$. Hence
\begin{eqnarray*} E_1(u) &=&
\eta\int_{-1}^1 |\ddot u |^2 dx - \int_{-1}^1 fu\, dx\,+\,\mu \int_{-1}^1|u-w|^2 dx \, +
J(u)
\,=
\\ &=&  J(u) +
\mu \int_{-1}^1 \!\! \big( (w-u)u + (u-w)^2 \big)dx 
+\,(-1\!+\!1/2)\!\int_{-1}^1 fu\, dx \,+
\\ &&
\,+\, \eta\,\Big(-\ddot u^+(-1) \dot u(-1) \,+\, \ddot u^-(1) \dot u(1)
+ \dddot u^+(-1) u(-1) \,-\, \dddot u^-(1) u(1)\Big)\,.
\end{eqnarray*}

\vskip0.3cm
Now assume $(u_r,u_p)\in\argmin F_1$ and label $t_{j}$ as above.
\\
By taking into account \eqref{EulerF1dprima}-\eqref{EulerF1endpointsSDumenoSu} and performing integrations by parts, we get
$$ \eta\!\int_{-1}^1\!|\ddot u_r|^2 dx\,=
\eta\!\int_{-1}^1\!\!(\ddot u_r)''u_r \,dx
\,=\, -\,\mu\!\int_{-1}^1\!\! (u_r-u_p)\,u_r \,dx
\,+ \int_{-1}^1 \!\!(f_r/2)\,u_r\,dx
\, +\, \eta \left[\ddot u_r \dot u_r\!-\!\dddot u_r\, u_r \right]_{-1}^{+1}  . $$
Performing two
integrations by parts and taking into account \eqref{EulerF1dseconda},
we get
\begin{eqnarray*}
     \gamma\! \int_{-1}^1\!|u_p''|^2 dx &=&
     \gamma\!\int_{-1}^1\!\!\!u_p''''u_p \,dx
     + \gamma\,\left[\,u_p''u_p'\!-u_p'''u_p\,\right]_{-1}^1
     \\
     &=&
\mu\!\int_{-1}^1\!\!(u_r-u_p)\,u_p\,dx\,  + \,\frac 1 2 \int_{-1}^1\!\!f_p u_p \,dx
+\gamma\,\!\left[\,u_p''u_p'\!-u_p'''u_p\,\right]_{-1}^{+1}
           \ .
\end{eqnarray*}
Then for any $(u_r,u_p)\in\argmin F_1$ we obtain
\begin{eqnarray*}
  F_1(u_r,u_p)
 \! \!&=&  \!\!
  J(u_r)\, +\,(-1\!+\!1/2)\!
  \left(
  \int_{-1}^1 f_ru_r\,dx\,+\, 
  \!\int_{-1}^1 f_pu_p\, dx \right)
  -\, \mu \int_{-1}^1  (u_r-u_p)^2  dx
  \\
   &&  + \, \left[ \,\eta \,(\,\ddot u_r \,\dot u_r- \dddot u_r\, u_r) +\gamma\,(\,u_p''w'-u_p'''w)\,\right]_{-1}^1
   \ . \quad \square
\end{eqnarray*}
\begin{prop}\phantom{.}\label{EulerThmConstr}
(\textbf{
Variational conditions for the minimizers of $E_1$ under unilateral constraint})
\\ 
Every minimizer $u$ of $E_1$ over the closed convex set
\begin{equation}
K := \left\{\,v\!\in\!  H^2(-2,2):\, v\!\ge\! w \,\hbox{ on} \, [-1,1], \,  v \!=\! w \, \hbox{ on}\, (-2,-1)\cup(1,2)\right\}
\end{equation}
(see Remark \ref{rmk constraint} about the pointwise everywhere meaning the unilateral  constraint) 
fulfils the variational inequality
\begin{equation}\label{EulerE1dConstr}
   \left\{
\begin{array}{l} 
\qquad u\in K: 
\\
\displaystyle
\int_{(-1,1)\setminus
     \left( S_{u}\cup S_{\dot u} \right)}
     \!\!
    \big(\eta \,u''''+\mu (u-w)-f/2\big)\, \big(u-v\big)
    \ \le\ 0  
    \qquad \forall v\in K\, ,
     \end{array}
     \right.
\end{equation}
together with the bilateral conditions at the free discontinuity and free-gradient discontinuity set where the contact does not play a role:
\begin{equation}\label{EulerE1dSuConstr}
    \ddot u^+ = \ddot u^-   = 0 \qquad\,
    \ \hbox{on }
      \big(S_{u}\!\setminus\!\{\pm 1\}\big) 
     \cap \{u^+ \!>\! w\}\cap \{u^-\!>\! w\},    
\end{equation}
\begin{equation}\label{EulerE1dSuConstrBis}
    \dddot u^+= \dddot u^-  = 0 \qquad\,
    \hbox{on }
     \big(S_{u}\!\setminus\!\{\pm 1\}\big) 
     \cap \{u^+ \!>\! w\}\cap \{u^-\!>\! w\},
\end{equation}
%
\begin{equation}\label{EulerE1dSDuConstrBis}
    \ddot u^+ = \ddot u^-  = [\dddot u]  = 0 \qquad \hbox{on }
     \big(S_{\dot u}\!\setminus\!\left( S_u\cup \{\pm 1\} \right) \big)
      \cap \{u^+\!>\!w\}\cap \{u^- \!>\!w\},
\end{equation}
\begin{equation}\label{EulerE1dsintesiConstr} 
\begin{array}{l}
\ddot u\in H^2\Big((-1,1)\setminus\{u^+=w \hbox{ or }u^-=w\}\Big)\quad\hbox{and}
\vspace{0.1cm}
\\
 \eta\,(\ddot u)''+\mu (u-w)=f/2 \quad \hbox{in }\
 \mathcal{D}'\big((-1,1)\!\setminus\! \{u^+=w \hbox{ or }u^-=w\}\big)\ ,
\end{array}
     \end{equation}
\begin{equation}\label{EulerE1endpointsSuConstr}
\left\{\
\begin{array}{l}
    \ddot u^+ (-1)= \dddot u^+(-1) =  0 \qquad \hbox{if }-1 \not\in S_{\dot u}\setminus S_u
    \hbox{ and } u^+(-1) \!>\! w(-1), \vspace{0.2cm} \\
    \ddot u^- (+1)= \dddot u^-(+1) =  0 \qquad \hbox{if }+1 \not\in S_{\dot u}\setminus S_u \hbox{ and } u^-(+1)\!>\! w(+1),
    \end{array} \right.
\end{equation}
\begin{equation}\label{EulerE1endpointsSDumenoSuConstr}
\left\{\
\begin{array}{l}
    \ddot u^+ (-1)=   0 \qquad \hbox{if }-1\in S_{\dot u}\!\setminus\! S_u \hbox{ and } u^+(-1) \!>\!w(-1), \vspace{0.2cm} \\
    \ddot u^- (+1)=   0 \qquad \hbox{if }+1\in S_{\dot u}\!\setminus\! S_u
    \hbox{ and } u^-(+1) \!>\! w(+1),
    \end{array} \right.\qquad\phantom{,}
   \end{equation}
jump condition $[\dddot u^+]=0$ in \eqref{EulerE1dSDuConstrBis} can be improved when $\alpha\!=\!\beta$, hence
\begin{equation}\label{EulerE1dSuimprovedConstr}
  \hbox{if }\alpha=\beta \hbox{ then } 
  \left\{ \
  \begin{array}{l}
  \dddot u^+ = \dddot u^- = 0 \quad
    \hbox{ on }\Big(\left( S_{u}\cup S_{\dot u} \right)\! \setminus\!\{\pm 1\} \Big)
     \cap \{u^+\!>\!w^+\}\cap \{u^- \!>\!w\},
    \\
    \dddot u^+(-1)  =  0  
    \quad\hbox{if }-1\in S_u \cup S_{\dot u}
    \hbox{ and } u^+(-1) \!>\! w(-1),
    \vspace{0.1cm}\\
    \dddot u^-(+1) =  0 
    \quad\hbox{if }+1\in S_u \cup S_{\dot u}
    \hbox{ and } u^-(+1) \!>\! w(+1),
    \end{array}
    \right.
\end{equation}
and in addition the unilateral conditions at the free discontinuity and free-gradient discontinuity sets of $u$ where the contact with the obstacle plays a role:
\begin{equation}\label{EulerE1dSuConstrBisContact+}
    \ddot u^+ \!\ge 0 \qquad\,
    \hbox{on }
     \big((S_{u}\cup S_{\dot u})\!\setminus\!\{ +1\}\big) 
     \cap \{u^+ \!=\! w\}\,,
     \end{equation}
\begin{equation}\label{EulerE1dSuConstrBisContact-}
    \ddot u^- \!\ge 0 \qquad\,
    \hbox{on }
     \big((S_{u}\cup S_{\dot u})\!\setminus\!\{- 1\}\big) 
     \cap \{u^- \!=\! w\}\,.
\end{equation}
No condition on $\dddot u^\pm$ is present on $S_{u}\cup S_{\dot u}$.
\end{prop}
\vskip0.1cm
\textsl{Proof - } The proof repeats the first 7 steps of Proposition \ref{EulerThm} proof, but achieves less information since there is a strictly smaller set of admissible variations.
\\
Step 1 is fully recovered thus, here we can exploit the Green formula \eqref{formula 2}.
\\ We repeat Steps 2-7, by performing all the admissible variations of $u$ which are of the kind $u+\varepsilon (v-u)$, with $\varepsilon\in[0,1]$ and $v\in K$: for comparison, here $\varphi=v-u$.
\\
As in Step 3 for the case of non constrained competitors, for
$\varphi =v-u$ belonging to $ \mathbf{H}^{2}\cap C^{2}([t_{l},t_{l+1}])$, $%
l=0,...,\mathsf{T}$, with $\mbox{\textnormal{spt}}(\varphi )\subset A$
Borel subset of $[-1,1]$ and $(S_{\dot{u}}\setminus S_{u})\cap A=\emptyset $
we still get
$\big(S_{u+\varepsilon \varphi }\cup S_{\dot{u}+\varepsilon \dot{\varphi}}\big)\cap A\ \subset \ {S_{u}}\cap A$
and
\begin{eqnarray}
\nonumber
0  \!\!& \leq &\!\!
   E_1(u+\varepsilon \varphi ,A)-E_1(u,A) \ =
\\
\nonumber
    \!\!& = & \!\! \alpha \left( \sharp (S_{u+\varepsilon \varphi }\cap A)-\sharp (S_{u}\cap
A)\right) +\beta \sharp \left( (S_{\dot{\varphi}}\setminus S_{u+\varepsilon
\varphi })\cap A\right) \,+
\\
\nonumber
    &  &  \qquad + \
   2\varepsilon \left( \sum_{l=0}^{\mathsf{T}}\int_{t_{l}}^{t_{l+1}}\big(\eta\, u^{\prime
\prime }\varphi ^{\prime \prime }+\mu(u-w)\varphi - (f/2) \varphi\big)\,dx\right) \, + \, o(\varepsilon )\, =
\\
\nonumber
    \!\!& = & \!\!
   \alpha \left( \sharp (S_{u+\varepsilon \varphi }\cap A)-\sharp (S_{u}\cap
A)\right) +\beta \sharp \left( (S_{\dot{\varphi}}\setminus S_{u+\varepsilon
\varphi })\cap A\right) \,+
\\
\nonumber
   &  & \qquad + \
   2\varepsilon \bigg(\sum_{l=0}^{\mathsf{T}}\int_{t_{l}}^{t_{l+1}}\big(\eta\, u^{\prime
\prime \prime \prime }\varphi +\mu(u-w)\varphi - (f/2) \varphi\big)\,dx \,+
\\
\nonumber
   &  & \qquad + \ {\dddot{u}}^{+}(-1)\varphi ^{+}(-1)-\ddot{u}^{+}(-1)\dot{\varphi}^{+}(-1)-{%
\dddot{u}}^{-}(1)\varphi ^{-}(1)+\ddot{u}^{-}(1)\dot{\varphi}^{-}(1) +
\\
\label{3.50}
   &  & \qquad + \ \,\eta\!\sum_{\left(S_{u}\cap A\right)\setminus\{\pm 1\}}\!\Big(\mathbf{\big[\kern-5pt\big[}\dddot{u}\varphi \mathbf{%
\big]\kern-5pt\big]}-\mathbf{\big[\kern-5pt\big[}\ddot{u}\dot{\varphi}%
\mathbf{\big]\kern-5pt\big]}\Big)\bigg)\ +\ o(\varepsilon )\,.
\end{eqnarray}
As in Step 4, let $\varphi \in \mathbf{H}^{2}\cap C^{2}([t_{l},t_{l+1}])$, $l=0,...,%
\mathsf{T}$, $\mbox{\textnormal{spt}}(\varphi )\subset A$, with $A$ Borel subset of $(-1,1) $ and $S_{\varphi
}=\emptyset =(S_{u}\setminus S_{\dot{u}})\cap A$. Then, up to a finite set of possible values of $\varepsilon $ entailing cancelation of $\dot u$ discontinuity, we can choose $\varepsilon $ arbitrarily
small such that
$
(S_{u+\varepsilon \varphi }\cup S_{\dot{u}+\varepsilon \dot{\varphi}})\cap
A=S_{\dot{u}+\varepsilon \dot{\varphi}}\cap A=S_{\dot{u}} \,
$; 
thus, by Green formula (\ref{formula 2})
\begin{eqnarray}
\nonumber
0 \!\!&\leq &\!\!
E_1(u+\varepsilon \varphi ,A)-E_1(u,A) \,\leq
\\
\nonumber
  \!\!& \leq &\!\! \beta \big(\sharp \,(S_{\dot{u}+\varepsilon \dot{\varphi}}\cap
A)-\sharp \ (S_{\dot{u}}\cap A)\big)\,+
\\
\nonumber
  && \qquad +\ 2\varepsilon\, \bigg(\sum_{l=0}^{\mathsf{T}}\int_{t_{l}}^{t_{l+1}}\big(\eta\, u^{\prime
\prime }\varphi ^{\prime \prime }+\mu(u-w)\varphi - (f/2) \varphi\big)\,dx\bigg)\,+\,o(\varepsilon ) \,=
\\
\nonumber
  \!\!& = &\!\!
  2\varepsilon \, \bigg(\sum_{l=0}^{\mathsf{T}}\int_{t_{l}}^{t_{l+1}}
  \Big( \eta\, u^{\prime
\prime \prime \prime }+\mu (u-w)- (f/2)\Big)\,\varphi \ dx\ +
\\
\nonumber
&& \qquad
+\ {\dddot{u}}^{+}(-1)\varphi ^{+}(-1)-\ddot{u}^{+}(-1)\dot{\varphi}^{+}(-1)-{%
\dddot{u}}^{-}(1)\varphi ^{-}(1)+\ddot{u}^{-}(1)\dot{\varphi}^{-}(1)\,+
\\
\label{3.51}
&& \qquad
+\ \eta \sum_{\left(S_{\dot{u}}\cap A\right)\setminus\{\pm 1\}}\Big(\mathbf{\big[\kern-5pt\big[}+\dddot{u}\varphi
\mathbf{\big]\kern-5pt\big]}-\mathbf{\big[\kern-5pt\big[}\ddot{u}\dot{\varphi%
}\mathbf{\big]\kern-5pt\big]}\Big)\bigg)\ +\ o(\varepsilon )\ .
\end{eqnarray}
In all cases now $\varphi=v-u$ with $v\in K$.
\\
By all choices of open sets $A$ and $v\in K$ fulfilling $\spt (v-u)\subset A \subset\subset (-1,1)\!\setminus\! (S_u\cup S_{\dot u})$
 we get,
\begin{eqnarray*}
&&\displaystyle
 \int_{(-1,1)\setminus
     \left( S_{u}\cup S_{\dot u} \right)}
     \!\! 
    \big(\eta u''''+\mu(u-w)-f/2\big)(v-u)\,dx
    =
    \\
    &&\displaystyle
    =\ \int_{(-1,1)\setminus 
     \left( S_{u}\cup S_{\dot u} \right)}
     \!\!
     \big(\eta u''(v-u)''+\mu(u-w)-f/2\big)(v-u)\,dx
     \ \ge\ 0  
\end{eqnarray*}
say \eqref{EulerE1dConstr}. Then, by inserting \eqref{EulerE1dConstr} in \eqref{3.50},\eqref{3.51}, we single out the conditions at every point of singular set.
\vskip0.1cm
\textbf{Outside the contact set} $\{u^+=w\}\!\cup\!\{u^-=w\}$ we can repeat the discussion made in the proof of Proposition \ref{EulerThm}, since 
$\varphi^\pm=(v-u)^\pm$ and $\dot\varphi^\pm=(\dot v-\dot u)^\pm$ are allowed to achieve both positive and negative values outside the contact set.
\\ 
Up to a finite set of possible values of $\varepsilon $ entailing cancellation of discontinuity, we have $S_{u+\varepsilon
\varphi }\cap A=S_{u}\cap A$. Then by discarding such values we can choose arbitrarily small $%
\varepsilon $ satisfying
\begin{equation*}
\sharp \,((S_{\dot{\varphi}}\setminus S_{u+\varepsilon \varphi })\cap
A)=\sharp \,((S_{\dot{\varphi}}\setminus S_{u})\cap A)=0
\end{equation*}
By taking into account \eqref{EulerE1d} and the arbitrariness of the two traces of $\varphi
$ and $\dot{\varphi}$ on the two sides of points in $S_{u}$, for small $%
\varepsilon $, we can choose $\varphi $ with $\varphi ^{\pm }=0$, and $\dot{%
\varphi}^{+}=0$ together with $\dot{\varphi}^{-}$ arbitrary, or viceversa to get $\,\ddot u^\pm=0$ on $S_u\setminus\{\pm 1\}$.\\ Similarly, we obtain $\,\dddot u^\pm=0$ on $S_u\setminus\{\pm 1\}$ by choosing $\dot{%
\varphi}^{\pm }=0$, and $\varphi ^{+}=0$ together with $\varphi ^{-}$ arbitrary or
vice-versa. So \eqref{EulerE1dSuConstr} is proved.\\
If some clamped endpoint ($-1$ and/or $+1$) belong to $S_u$, then \eqref{EulerE1endpointsSuConstr} is obtained as above, but taking into account that $\varphi\equiv 0$ outside $[-1,1]$.\\
Summarizing, we obtain \eqref{EulerE1dSuConstrBis},
\eqref{EulerE1dsintesiConstr},\eqref{EulerE1endpointsSuConstr},\eqref{EulerE1endpointsSDumenoSuConstr}, hence \eqref{EulerE1dSDuConstrBis} and \eqref{EulerE1dSuimprovedConstr}, by the same argument of Steps 3\,-\,7.
\vskip0.2cm
\textbf{On the contact set} $\{u=w\}$ we can repeat again the discussion made in the proof of Proposition \ref{EulerThm}, but here the coefficient of $2\varepsilon$ in \eqref{3.51} must be only nonnegative, thus we get inequalities in place of equalities. Moreover, since 
$\dot \varphi^+=(\dot v-\dot u)^+$ is allowed to achieve only positive values
 and 
$\dot \varphi^-=(\dot v-\dot u)^-$ is allowed to achieve only negative values, whereas left and right values have always opposite sign,
we deduce
\eqref{EulerE1dSuConstrBisContact+}, \eqref{EulerE1dSuConstrBisContact-}. On the other hand, on the contact set
$\varphi^\pm=( v- u)^\pm$ is always null; therefore, we get no condition on every term whose multiplier is $\varphi^+$ or $\varphi^-$. \hfill $\square$
\begin{prop}\phantom{.}\label{EulerThmConstrF1}
(\textbf{
Variational conditions for the minimizers of $F_1$ under unilateral constraint})
\\ 
Every minimizing pair $(u_r,u_p)$ of $F_1$ over the convex set
\begin{equation*}
\mathbf K \!:=\! \left\{(v_r,v_p)\!\in\! \mathbf  H ^2(-2,2)\!\times\! H^2(-2,2)\!:\, v_r\!\ge\! v_p \,\hbox{on} \, [-1,1], \,  v_r \!=\!v_p\!=\! w \, \,\hbox{on}\, (-2,-1)\!\cup\!(1,2)\right\}
\end{equation*}
fulfils the quasi-variational inequalities
\begin{equation}\label{EulerF1dConstr ur}
   \left\{
\begin{array}{l} 
 u_r\in\mathbf H^2(-2,2): \ v_r\!\ge\! v_p \,\hbox{on} \, [-1,1], \,  v_r \!=\!v_p\!=\! w \, \,\hbox{on}\, (-2,-1)\!\cup\!(1,2) \hbox{ and}
 \vspace{0.1cm}
\\
\displaystyle
\!\!\int_{(-1,1)\setminus
     \left( S_{u_r}\cup S_{\dot u_r} \right)}
     \!\! \!\! \!\! \!
    \big(\eta \,u''''_r\!+\mu (u_r-u_p)-f_r/2\big)\, \big(u_r-v\big)
    \, \le\, 0  
    \quad \forall v: (u_r,v)\!\in\! \mathbf K ,
     \end{array}
     \right.
\end{equation}
\begin{equation}\label{EulerF1dConstr up}
   \left\{
\begin{array}{l} 
 u_p\in H^2(-2,2): \ v_p\!\le\! v_r \,\hbox{on} \, [-1,1], \,  v_r \!=\!v_p\!=\! w \, \,\hbox{on}\, (-2,-1)\!\cup\!(1,2) \hbox{ and}
 \vspace{0.1cm}
\\
\displaystyle
\!\!\int_{(-1,1)\setminus
     \left( S_{u_r}\cup S_{\dot u_r} \right)}
     \!\! \!\! \!\! \! \!
    \big(\gamma \,u''''_p\!+\mu (u_p-u_r)-f_p/2\big)\, \big(u_p-z\big)
    \, \le\, 0  
    \quad \forall z: (z,u_p)\!\in\! \mathbf K ,
     \end{array}
     \right.
\end{equation}
together with the standard bilateral conditions (say \eqref{EulerE1dSuConstr}-\eqref{EulerE1endpointsSDumenoSuConstr}
 with $u_r,u_p$ replacing respectively $u,w$) at the free discontinuity and free-gradient discontinuity set where the contact does not play a role, 
and the unilateral conditions at the free discontinuity and free-gradient discontinuity set where the contact with the obstacle plays a role:
\begin{equation}\label{EulerE1dSuConstrBisContact+ F1}
    \ddot u_r^+ \!\ge 0\qquad 
    \hbox{on }
     \big((S_{u_r}\cup S_{\dot u_r})\!\setminus\!\{ +1\}\big) 
     \cap \{u_r^+ \!=\! u_p\}\,,
     \end{equation}
\begin{equation}\label{EulerE1dSuConstrBisContact- F1}
     \ddot u_r^- \!\ge 0\qquad 
    \hbox{on }
     \big((S_{u_r}\cup S_{\dot u_r})\!\setminus\!\{- 1\}\big) 
     \cap \{u_r^- \!=\! u_p\}\,.
\end{equation}
No condition on $\dddot u_r^\pm$ is present on $S_{u_r}\cup S_{\dot u_r}$.
\end{prop}
\textsl{Proof - }
Repetition of the steps of the last proof provides the proof the claims about minimizers $(u_r,u_p)$ of $F_1$ with unilateral implicit constraint,
by performing all the admissible variations of $u_r$ which are of the kind 
$u_r+\varepsilon (v-u_r)$, with $\varepsilon\in[0,1]$ 
and $(v,u_p)\in \mathbf K$ and
$v_r+\varepsilon (z-v_r)$, with $\varepsilon\in[0,1]$ and and $(u_r,z)\in \mathbf K$. \hfill $\square$
\section{Hard-device reinforcement of flexural plate}\label{hard-device Section}
In this section we deduce the existence statement in the case of hard-device reinforcement: minimization of functional $E$ defined by \eqref{E}.
\vskip0.1cm
\textsl{Proof of Theorem 2.1} -
After noticing that by
\begin{equation*} 
\, E(\emptyset,w)\,=\, \eta\|D^2w\|_{L^2(\Omega)}^2 - \int_\Omega fw \,d\xx\,<\,+\infty\ .
\end{equation*}
the domain of $E$ is not empty, and by
\begin{equation}
\label{reduction E to BZ}
    \int_\Omega \!\!\big( \,\mu\,  |v-w|^2 \, -f v\,\big)\, dx =
     \int_\Omega \!\! \mu\Big(  v-\big(w+f/(2\mu)\big)\Big)^2dx
     \,-
     \int_\Omega \!\!\Big( fw+f^2/(4\mu)\Big)dx
               \end{equation}
where the last summand on the right-hand side is a constant,
we have that the functional $E$ is bounded from below since
beside such constant all other terms are nonnegative.

The notion of essential admissible pairs, set by \eqref{admissible K v}, selects (\cite{CLTMJM})
those pairs $(K,v)$ which are cleansed of every spurious artifact
that does not affect the functional value and are good representatives in equivalence
classes of admissible pairs. This definition of admissible pair prevents diffused damage but
allows to prove partial regularity of displacements $v$:
free discontinuity (crack) and free gradient discontinuity (folds)
are allowed in competing configurations of the structure.
\\
Thus, the claims of present Theorem \ref{Theorem 1} follow from Theorem 2.3
in \cite{SBZ} and \cite{SBZDIR} about functional (2.2) defined therein,
by setting $g=w+f/(2\mu)$, a datum which belongs to $ L^4(\Omega)$ due to present assumptions. Precisely
we can choose
$\widetilde\Omega=\Omega_p$, $\alpha=\beta=1$; hence (2.3),(2.4),(2.5) and (2.20) of \cite{SBZDIR} are fulfilled
thanks to the conditions \eqref{Omega},\eqref{eta mu} assumed here.
Moreover $D^2w\in L^{\infty}(A)$ for any open set s.t.
$\Omega\subset\!\subset A \subset\!\subset \Omega_p$ and we have that the sets $M,T_0,T_1$ (as denoted in \cite{SBZDIR})
are empty, hence
(2.6)-(2.11) of \cite{SBZDIR} hold true thanks to the assumption \eqref{w} made here.

\hfill $\square$

\section{Strengthening reinforcement of flexural plate}\label{strengthening Section}
In this section we deduce the existence statement in the case of strengthening reinforcement: minimization of functional $F$ defined by \eqref{F}.
\\
To deal with the case of strengthening reinforcement we need a relaxed formulation of functional \eqref{F}, as it is usual in the analysis of free discontinuity problems. 
We list standard notations 
(see \cite{AFP},\cite{SBZ},\cite{BZEE},\cite{CLTadvcv}):
$B_\varrho (\mathbf{x})$ denotes the open ball
$\{\,\mathbf{y}\!\in\!\R^2 \!:\, |\mathbf{y}-\mathbf{x}| \!<\! \varrho\,\}$;

$\mathcal{H}^1(A)$ and 
$|A|$ denote respectively, the 1-dimensional Hausdorff measure
and the outer Lebesgue measure of a subset $A\subset\R^2$;
for every Borel function $v:\Omega\to\R$ and $\mathbf{x}
\in \Omega$, $z \in \overline\R :=  \R \cup \{-\infty,+\infty\}$,
we set $\disp{z} = {\rm ap}\lim_{\mathbf{y}\to \mathbf{x}}
v(\mathbf{y})\,$  
\vskip-0.25cm 
(notation for the 
{\sl approximate limit} of $v$ at $\mathbf{x}$)
if,
for every $g\in C^0 (\overline\R),$
$$
   g({z}) \ =\  \lim_{\ro \to 0}\, 
    |B_\ro(\mathbf{0})|^{-1} \!\!\int_{B_\ro(\mathbf{0})}
   \!\! g(v(\mathbf{x}+\mbox{\boldmath$\xi$}))\,d\mbox{\boldmath$\xi$}\, ;
$$
the function $\disp\widetilde v(\mathbf{x})=\mathop{\rm ap\,lim}\limits_{\mathbf{y}\to
\mathbf{x}} v(\mathbf{y})$ is called {\sl representative} of
$v\,;$
$$ \disp S_v =
\{\mathbf{x}\in\Omega:\,\not\!\exists z\ \mbox{ such that }\
\mathop{\rm ap\,lim}\limits_{\mathbf{y}\to \mathbf{x}} v(\mathbf{\mathbf{y}})=z\}
\quad\hbox{is the {\sl singular set} of }v\,.$$
A Borel function $v:\!\Omega\!\to\!\R$ is {\sl approximately continuous at} $\mathbf{x}
\!\in\! \Omega$ iff $\disp v({\mathbf x})\!=\!\mathop{\rm ap\,lim}\limits_{\mathbf{y}\to
\mathbf{x}} v(\mathbf{y})$.
\\
As usual, $Dv$ denotes the
distributional gradient of $v$ and $\nabla v(\mathbf{x})$ denotes
the \textsl{approximate gradient} of $v\,,$ say
$v$ is approximately differentiable at $x$ if there
exists a vector $\nabla v(\xx)\in\R^2$ (the approximate gradient of
$v$ at $\xx$) such that
\[\mathop{\rm ap\,lim}\limits_{\mathbf{y}\to \mathbf{x}}\frac{|v(\mathbf{y})-\widetilde{v}(\mathbf{x}) - \nabla v(\mathbf{x})\cdot
(\mathbf{y}-\mathbf{x})|}{|\mathbf{y}-\mathbf{x}|}=0.\]
A function $u\!\in\! BV(\O)$ is approximately differentiable a.e., moreover 
for $\h$ almost every $\mathbf{x}\!\in\! S_u$ there exist
$\nu(\mathbf{x})\!\in\partial \!B_1$, $v_+(\mathbf{x})\!\in\!\R$, $v_-(\mathbf{x})\!\in\!\R$ with
$v_+(\mathbf{x})\!>\!v_-(\mathbf{x})$ such that 
\[\lim_{\varrho \to 0} \ \varrho ^{-n}\!\int_{\{\mathbf{y}\in B_{\varrho };\, \mathbf{y}\cdot \nu(\mathbf{x})>0\}}
|v(\mathbf{x}+\mathbf{y})-v_+(\mathbf{x})|\,d\mathbf{y}=0,\]
\[\lim_{\varrho \to 0} \ \varrho ^{-n} \!\int_{\{\mathbf{y}\in B_{\varrho };\, \mathbf{y}\cdot \nu(\mathbf{x})<0\}}
|v(\mathbf{x}+\mathbf{y})-v_-(\mathbf{x})|\,d\mathbf{y}=0.\]
$\sbvo$ denotes the De Giorgi class of functions $v\in BV (\Omega)$ such that
$$
\int_\Omega|Dv| = \int_\Omega|\nabla v| \, d\mathbf{x} + \int_{S_v}| v^+ -
v^-| \, d\h .$$
We introduce:
$$\sbv_{loc}(\Omega)\ :=\ \left\{v\in
\sbv(\Omega');\  \forall\,\Omega'\subset\kern-1pt\subset
\Omega\right\},$$
\begin{equation}\label{gsbv}
\gsbv(\Omega) \ :=\ \big\{ v :{} \Omega\to\R \ \hbox{Borel
function}; -k \lor v \land k \in\sbv_{loc}(\Omega)\ \forall k\in
\mathbb{N} \big\}.
\end{equation}
\begin{equation}\label{gsbv2}
\gsbv^2(\Omega)\ :=\  \big\{\,v\in\gsbv(\Omega) ,~\nabla v\in \big(\gsbv(\Omega)\big)^2 \big\}.
\end{equation}
If $v\in \gsbv(\Omega)$ then $\nabla v$ exists a.e., and for $v\in \gsbv^2(\Omega)$ we set
$\nabla^2 v=\nabla (\nabla v)$. 

\vskip0.4cm
Eventually, we introduce the weak formulation $\mathcal{F}$ of functional $F$ defined by \eqref{F}:
\begin{equation}\label{weak F}\left\{
\begin{array}{l}     \mathcal{F}(v_r,v_p) \, =\, \mathcal{F}(V)
    \ :=\   \\  \displaystyle \ \
    \mathcal{H}^1\!\left( S_{v_r}\right) +\,
               \eta\! \int_{\Omega}\! |\nabla^2 v_r|^2  d\xx +
               \mu\!\int_{\Omega}\! |v_r-v_p|^q \, d\xx +
               \!\int_{\Omega} \!\!\!\left(\gamma\,|D^2 v_p|^2
               -  f \, v_p \right)  d\xx \ ,
               \vspace{0.1cm}
                \\ 
    \  \ \forall\, V=(v_r,v_p)\in \mathbb{X}:=
    \Big(GSBV^2(\Omega_p)\cap L^2(\Omega_p)\Big)\!\times \!H^2(\Omega_p)\,
    \\ \qquad\qquad\qquad\qquad\qquad\qquad\qquad
    \hbox{s.t.}\ \ v_r=v_p=w\ \hbox{ a.e. } \Omega_p\setminus \overline\Omega\,.\ \,
   \end{array}\right.
   \end{equation}
We emphasize that, since $v_p=w$ in $\Omega_p\setminus\overline\Omega$ and $w\in C^2(\Omega_p\setminus\overline\Omega),$ we get $S_{v_p}\cup S_{\nabla v_p}=\emptyset$ and
\begin{equation}\label{cal F + M + Fp}
\mathcal{F}(v_r,v_p) \ = \ \mathcal{F}_r(v_r) + M(v_r-v_p) + F_p(v_p)
\end{equation}
where $M$, $F_p$ and $\mathcal{F}_r$ are defined by \eqref{MM}, \eqref{Fp} and
\begin{equation}\label{weakFc}
    \mathcal{F}_r(v_r)\ :=\ \mathcal{H}^1\left( S_{v_r} \cup S_{\nabla v_r} \right)\, +\,
                \int_{\Omega}  \Big( 
                \eta|\nabla^2 v_r|^2 - f_rv_r
                \Big) \,d\xx \ .
\end{equation}
\begin{thm}\label{thm weak F}
Assume \eqref{Omega},\eqref{w},\eqref{eta mu} and \eqref{f}.
\\ 
Then
the functional $\mathcal{F}$ achieves a finite minimum over $\mathbb{X}\,.$
\end{thm}
\textsl{Proof -}
First we notice that $\mathcal{F}$ has non empty domain: in fact
\eqref{w} entails $S_w=S_{\nabla w}=\emptyset$ and
\begin{equation}\label{dom cal F not empty}
\mathcal{F}(w,w)\,<\, (\eta+\gamma)\|D^2w\|_{L^2(\Omega)}^2 - \int_\Omega (f_r+f_p)\,w \ d\xx\,<\,+\infty\ .
\end{equation}
We have the identity
\begin{equation}
\label{reduction F to BZ}
    \int_\Omega \!\!\big( \,\mu\,  |v_r-v_p|^2  -f_r v\,\big) dx =\!
     \int_\Omega \!\! \mu\Big(  v_r-\big(v_p+f_r/(2\mu)\big)\Big)^2dx
     - \!
     \int_\Omega \!\!\Big( f_r v_p+f_r^2/(4\mu)\Big)dx .
               \end{equation}
\\
If $K_\Omega$ denotes the best Poincar\'e inequality constant in $H^2_0(\Omega)$, namely
\begin{equation}\label{Poincare}
\|v\|_{L^2(\Omega)}^2 \,\le\, K_\Omega \,\|v''\|_{L^2(\Omega)}^2 
 \qquad \forall\,
v\!\in\! H^2_0(\Omega)\,,
\end{equation}
and, arguing as like as in \eqref{F1 inf boundedness}
we get, for everyy $V=(v_r,v_p)\in \mathbb X$
\begin{equation}
\label{F inf boundedness}
F(v_r,v_p) \ge 
        J(v_r) +
                \int_{-1}^1\! \!
                \Big(\eta|\ddot v_r|^2 +
               (\gamma/2)\,|v_p''|^2 \Big) \, dx \,
                    - C(\mu,\gamma,f_r,f_p,K_\Omega)\ .
\end{equation}
Then the functional $F$ is bounded from below on its domain.
Hence we can select a minimizing sequence $V_h=\big((v_r)_h,(v_p)_h\big)$ for $\mathcal{F}$\,: $\lim_h \mathcal{F}(V_h)=\inf \mathcal{F} \in \R$\,. \\ 
Thanks to \eqref{dom cal F not empty},
we may suppose
that
\begin{equation}\label{stima a priori su F}
c\,\leq\,\mathcal{F}(V_h)\leq C := \mathcal{F}(w,w)<+\infty \,.\end{equation}
Summarizing 
$F_p\big((v_p)_h\big)\leq C$, $F_r\big((v_r)_h\big)\leq C$
and  
$(v_p)_h$ is bounded in $H^2(\Omega)$.
Moreover there is $u_p\in H^2(\Omega_p)$ such that, up to subsequences and without relabelling,
$(v_p)_h$ is converging to $u_p$ weakly in $H^2(\Omega)$ and strongly in $L^2(\Omega)$, 
and $\|(v_p)''_h\|_{L^2} \to \|u_p''\|_{L^2}$\,.
\\
By using any fixed $(v_p)_h$ chosen from the sequence (which is bounded in $H^2$) as datum we find a minimizer, denoted by $z_h$, in
$GSBV^2(\Omega_p)\cap L^q(\Omega_p)$
of
\begin{equation*}
   v \ \mapsto\  \mathcal{F}_r(v) + M(v-(v_p)_h)
\end{equation*}
since this problem is equivalent to the minimization of Blake \& Zisserman functional for image segmentation with gray-level
datum $g=(v_p)_h+f_r/(2\mu)$
and Dirichlet boundary condition, referring to notation of Theorem 3.1 in \cite{SBZDIR}. Then
\begin{equation}\label{stima a priori su Fc}
 \mathcal{F}_r\big(z_h\big) + M\big(z_h-(v_p)_h\big) \ \leq 
 \mathcal{F}_r\big( (v_r)_h\big) + M((v_r)_h-(v_p)_h)
 \qquad \forall h\,.
 \end{equation}
Hence, by \eqref{cal F + M + Fp} and standard lower semicontinuity of $F_p$, the sequence of pairs $(z_h,(v_p)_h)$ is a minimizing sequence for $\mathcal{F}$ too.
Moreover, by \eqref{stima a priori su F},\eqref{stima a priori su Fc} we get
$$\mathcal{F}_r\big(z_h\big) + M\big(z_h-(u_p)_h\big) \leq C
$$
By compactness property of Theorem 8 in \cite{WBZ}, there are
$u_r\!\in\! GSBV^2(\Omega_p)\cap L^2(\Omega_p)$ and a subsequence s.t., again by extracting without relabeling,
$z_h\to u_r$. Moreover by lower semi-continuity property of Theorem 10 in \cite{WBZ}, we get:
$\,z_h\to u_r$ a.e., $z_h\rightharpoonup u_r$ $L^2$ and strongly in $L^s,$ $1\!\leq\! s <2,$
$\nabla z_h\to \nabla u_r$ a.e., $\nabla^2 z_h\to \nabla^2 u_r$ a.e., $\nabla^2 z_h\to \nabla^2 z$
weakly in $L^2$ and
\begin{equation*}
\mathcal{F}_r(u_r) + M(u_r-{u_p}) \leq \liminf_h \mathcal{F}_r(z_h) + M\big(z_h-(v_p)_h\big)
 \leq  
 \, C \ .
\end{equation*}
\vskip-0.2cm
Thus the pair $(u_r,u_p)$ is a minimizer of relaxed functional $\mathcal F$.
\hfill $\square$
\vskip0.3cm
\textsl{Proof of Theorem 2.2} -
Let $V=(u_r,u_p)\in\argmin \mathcal{F}$ (the existence of at least one such $V$ is warranted by Theorem \ref{thm weak F}). Then $u_p$ minimizes $z\mapsto F_p(z)+ M(v_r-z)$ among
$z\in H^2(\Omega_p)$ s.t. $z=w$ a.e. $\Omega_p\setminus\overline\Omega $.
So, due to \eqref{w} and \eqref{Omega}, $u_p\in C^2\cap L^\infty(\Omega)$. \\
Moreover, if $V=(u_r,u_p)\in\argmin \mathcal{F}$, 
then $u_p\in H^4(\Omega)$ and, referring to \eqref{MM} and \eqref{Fp}, $u_r$ minimizes $v\mapsto \mathcal{F}_r(v)+ M(v-u_p)$ among $v\in \gsbv^2(\Omega_p)\cap L^2(\Omega_p))$ s.t. $v=u_p=w$ a.e. $\Omega_p\setminus\overline\Omega $.
Thus, 
exploiting the identity \eqref{reduction F to BZ},
by Theorem 2.2 of \cite{SBZDIR} with the choices $\alpha=\beta=1$, $g=u_p+f_r/(2\mu)\in  L^4(\Omega_p)$ and $M=T_0=T_1=\emptyset$, and setting $Z=\overline{S_{u_r}\cup S_{\nabla u_r}}$, we obtain that the triplet $(Z,\widetilde{u_r},u_p)$ is an essential admissible triplet that minimizes $F$. 
\\ 

By applying the regularization argument detailed in \cite{SBZ},\cite{SBZDIR},\cite{CLTadvcv}
we obtain that
$\widetilde{u_r}\in C^2(\Omega_p\setminus Z)$, where
$Z$ is the smallest closed subset of $\Omega_p$ containing the region where $C^2$ regularity of $ \widetilde{u_r}$ is missing, and
$\mathcal{H}^1(Z\setminus (S_{u_r}\cup S_{\nabla u_r})=0$.
%
Eventually
\begin{equation*}
\begin{array}{l}
\hskip-0.2cm\mathcal{F}_r(\widetilde{u_r}) + G(\widetilde{u_r}-u_p) + F_p(u_p)
\leq \mathcal{F}_r({u_r}) + G({u_r}-u_p) + F_p(u_p)\leq
\\
\qquad\qquad\qquad\qquad\qquad\qquad\quad\, 
\leq \liminf_h \left(\
\mathcal{F}_r(z_h) + G(z_h-u_h) + F_p(u_h)\ \right)= \inf_{\mathbb X} \mathcal{F}
\end{array}
\end{equation*}
hence $\mathcal{F}(\,\widetilde{u_r}\,,\,u_p\,)\,=\,\min_{\mathbb X }\, \mathcal{F}\,.$
\\

Summarizing 
$F\big(\,Z,\,\widetilde{u_r},\,u_p\,\big)=\min \{\, F(K,v_r,v_p): (K,v_r,v_p) \hbox{ admissible triplet} \,\}\,.$
\hfill $\square$
\vskip0.3cm
\begin{rmk}\label{rmk 5 omega}
We emphasize that also the non-interpenetration between plate and reinforcement could be taken into account: e.g., adding the constraint $v\geq w  \ a.e\ \Omega_p$ to the essential admissible pairs for hard-device reinforcement 
and adding the constraint $v_r\geq v_p \ a.e\ \Omega_p$ to the essential admissible triplets for strengthening reinforcement.
Notice that here Remark \ref{rmk constraint} does not apply: competing functions are functions defined only almost everywhere, therefore the unilateral constraints act in the almost everywhere sense only.
These unilateral constraints do not introduce any additional difficulty in the study of the weak
formulations of both $F$ and $E$, since inequalities are preserved by compactness properties of minimizing sequences. 
Therefore Theorem \ref{thm weak F} holds true also under the additional constraint $v_r\ge v_p$.
\\ 
But the subsequent step required to show Theorem 2.2, say the proof of partial regularity for weak minimizers, would be not straightforward.
\\
For this reason in this short note we skip this substantial difficulty, postponing the analysis of the 2\,dimensional problems with unilateral constraints to a forthcoming paper.\\
However, in 1 dimension the strong and weak formulation do coincide,
so the analogous of Theorems \ref{Theorem 1} and \ref{Theorem 2} hold true
with or without the non-interpenetration constraint for beams: we have taken into account these constraints in the one-dimensional case by Theorems \ref{Thm 3.5}, \ref{Thm reinf1 with unilat constr} and Propositions \ref{EulerThmConstr}, \ref{EulerThmConstrF1}.
\end{rmk}
\section{Elastic-plastic reinforcement of flexural plate}\label{elastic-plastic Section}
In this section we deduce the existence statement in the case of strengthening reinforcement: minimization of functional $G$ defined by \eqref{G}.
\vskip0.2cm
\textsl{Proof of Theorem 2.3} -
By $G(\emptyset, w,w)\!=\! 
(\eta+\gamma) \|D^2 w\|^2_{L^2(\Omega)} \!- \int_\Omega (f_r+f_p) w \,d\xx 
< +\infty$, we know that the functional $G$ has nonempty domain.
Moreover
\eqref{F inf boundedness} warrants that the functional $G$ is bounded from below. 
\\
The existence of a minimizer of $G$ over 
over \textsl{essential admissible triplets} $(K,v_r,v_p)$, namely triplets fulfilling \eqref{admissible K vr vp G},
can be achieved by repetition of the direct method approach with the techniques of \cite{SPLATE}.
\\
Actually here, about minimization with respect to $v_r$, we have these differences with respect to \cite{SPLATE}: 
presence of the additional coupling term $\mu\int_\Omega |v_r-v_p|^2dx$; there are neither vanishing moments nor a safe load condition for the load $f_r$; last, there is a Dirichlet datum $w$ at the boundary. 
\\
However vanishing moments and load $f$ were exploited in \cite {SPLATE} only to achieve 
the boundedness from below of the functional, whereas here such boundedness is already
warranted by \eqref{F inf boundedness}. Moreover the additional term is a lower order perturbation, not affecting the existence of weak minimizers (thanks to the identity \eqref{reduction F to BZ}, still valid in present case), but requiring a technical correction in the proof of strong solutions by regularization of weak solutions.
\\
Precisely, first step (existence of weak solutions) requires no change: we introduce the space $SBH$ of Special Bounded hessian functions
$$SBH(\widetilde\Omega):=\{v\in H^{1,1}(\widetilde\Omega): Dv\in SBV(\Omega),\, v=w\hbox{ on }\widetilde\Omega\!\setminus\!\Omega)\}\,,$$ here  
$SBV$ is the space of bounded variation functions whose derivative has no Cantor part (\cite{AFP}); then we set the weak formulation of functional $G$ defined in \eqref{G}, defined on $v\in SBH(\Omega)$:
\vskip-0.7cm
\begin{eqnarray}\label{mathcal G}
     &\mathcal G(v_r,v_p) \!\!&:=\,
    \mathcal{H}^1\left( S_{Dv_r}\right)
     \,+\,
               \sigma \int_{ S_{Dv_r}} \!\! |  [Dv_r]\,|\,d\mathcal H^1
    \, +
                \int_{\Omega}\! \big(\eta |\nabla^2 v_r|^2 - f_r v_r\big)\, d\xx 
              +
              \\ && \nonumber
               \qquad \,+\,
               \mu\int_{\Omega} |v_r-v_p|^2 \, d\xx
                +\,
                \int_{\Omega}\! \big(\gamma |D^2 v_p|^2 - f_p v_p\big)\, d\xx
   \end{eqnarray}
\vskip-0.3cm
where $[z]$ denotes the jump of $z$
and 
$\nabla z$ denotes the approximate gradient of $z$, say the absolutely continuous part of $Dz$. 
\\
The existence of a minimizing pair $(u_r,u_p)$ for $\mathcal G$
follows by the same argument of present Theorem \ref{thm weak F}, taking into account of Theorem 2.9 in \cite{SPLATE}.
\vskip0.2cm
The proof of partial regularity in $\Omega$ for weak minimizers is achieved by exploiting blow-up and
quasi-minimizers as in Theorem 4.15 in \cite{SPLATE}: only Lemma 4.3 of \cite{SPLATE} must be adapted as detailed below, to take into account of the additional glueing term.\\
Still by identity \eqref{reduction F to BZ}, the load and glue terms together are represented  (up to the 
\vskip-0.1cm
addition of a constant irrelevant in minimization) by 
$ \displaystyle\mu \int_\Omega \! (  v-h)^2dx $ where $h:=u_p+f_r/(2\mu)$
and $h\in L^s, s>2$: this contribution replaces here the term $-\int_\Omega gv$ of \cite{SPLATE}, however this does not affect regularization of weak solutions, since, setting $\mathcal E(v)=\mathcal G(v)-\mu\int_\Omega (v-h)^2,$
every local minimizer of $\mathcal G$ is a local quasi minimizer of $\mathcal E$, due to the excess estimate (consequence of $u,v\in SBH(\Omega) \subset L^\infty(\Omega) $, $h\in L^s(\Omega)$, $s>2$ and H\"older inequality): 
\begin{equation*}
 \int_{B_\varrho(\xx)} \Big( (v-h)^2-(u-h)^2 \Big)\,dx\,=\,
  \int_{B_\varrho(\xx)} \Big( u^2-v^2-h(v-u) \Big) \, \le\, C\,\varrho^{2-2/s}
\end{equation*}
valid for $\overline{B_\varrho(\xx)}\subset\Omega$, $0<\varrho<1$ and $u,v\in SBH(\Omega)$ s.t. $v=u$ on $\Omega\!\setminus\!B_\varrho(\xx)$.
\\
Partial regularity at the boundary under Dirichlet condition, can be achieved by the same argument of \cite{SBZDIR}, taking into account of the simplifications due to the fact that here the competing functions are not only in $GSBV^2(\Omega)$, but they belong to $SBH(\Omega)$, hence they are globally continuous.
\\
Summarizing a minimizing pair $(u_r,u_p)$ of $\mathcal G$ leads to an essential minimizing triplet $\big( \overline{S_{\widetilde{Du_r}}}, u_r,u_p  \big)$ of $G$.
\hfill $\square$

\end{document}